%% file: minimax_expmix.tex
\documentclass[a4paper,11pt]{amsart}

\usepackage[english]{babel} 
\usepackage[latin1]{inputenc} 
\usepackage[T1]{fontenc} 
\usepackage{graphics}
\usepackage{times} 
\usepackage{amsfonts} 
\usepackage{multicol} 
\usepackage{amssymb} 
\usepackage{amsmath}
\usepackage{bbm} 
\usepackage{enumerate}
\usepackage{color, graphicx}
\usepackage{hyperref}%,showlabels}
\usepackage[numbers]{natbib}

\newtheorem{lem}{Lemma} 
\newtheorem{prop}{Proposition} 
\newtheorem{cor}{Corollary}
\newtheorem{theo}{Theorem}
\theoremstyle{definition}
\newtheorem{definition}{Definition} 
 
\theoremstyle{remark}
\newtheorem{remark}{Remark}

\newtheorem{ass}{Assumption}

\def\P{\mathbb{P}} 
\def\Xset{\mathrm{X}}
\def\Xfield{\mathcal{X}}
\def\E{\mathbb{E}} 
\def\N{\mathbb{N}}
\def\R{\mathbb{R}} 
\def\H{\mathbb{H}}
\def\1{\mathbbm{1}} 
\def\rme{\mathrm{e}}
\def\rmi{\mathrm{i}}
\def\rmd{\mathrm{d}}

\def\coeff{c}

\title[Nonparametric estimation of the Mixing Density]
{Nonparametric estimation of the Mixing Density using polynomials}

\author{Tabea Rebafka, Fran\c cois Roueff}

\begin{document}

\maketitle
% indicate corresponding author with \corref{}
% \author{\fnms{John} \snm{Smith}\corref{}\ead[label=e1]{smith@foo.com}\thanksref{t1}}
% \thankstext{t1}{Thanks to somebody} 
% \address{line 1\\ line 2\\ printead{e1}}
% \affiliation{Some University}

\section*{Abstract} 
We consider the problem of estimating the mixing density $f$ from $n$
i.i.d. observations distributed according to a mixture density with unknown
mixing distribution. In contrast with finite mixtures models, here the
distribution of the hidden variable is not bounded to a finite set but is
spread out over a given interval. 
We propose an approach to construct an  orthogonal series estimator of the mixing
density $f$  involving Legendre polynomials. 
The construction of the orthonormal sequence
varies from one mixture model to another. 
Minimax upper and lower bounds of the
mean integrated squared error are provided which apply in various contexts. 
In the specific case of exponential mixtures, it is shown that the 
estimator is adaptive over  a collection of specific smoothness
classes, more precisely,  there
exists a constant $A>0$ such that, when the order $m$ of the projection estimator verifies   $m\sim A \log(n)$, the 
estimator achieves the minimax rate   over this collection.
Other cases are 
investigated such as Gamma shape mixtures and scale mixtures of compactly
supported densities including Beta mixtures.
Finally,  a consistent estimator of the support of the mixing density $f$ is provided.

\section{Mixture Distributions}
We consider mixture distributions of densities belonging to some parametric
collection $\{\pi_t,t\in\Theta\}$ of densities with respect to the dominating
measure $\zeta$ on the observation space $(\Xset,\Xfield)$. A general representation of a mixture
density uses the so-called \textit{mixing distribution} and is of the following
form
\begin{equation}\label{def_mixture_dens}
\pi_{f}(x) = \int_\Theta f(t)\pi_t(x)\mu(\rmd t)\;,
\end{equation}
where the \textit{mixing density} $f$ is a density with respect to some measure $\mu$ defined on $\Theta$.  If $\mu$ is a
counting measure with a finite number of support points $\theta_k$, then obviously, $\pi_f$ is a finite mixture distribution
of the form $\sum_{k=1}^Kp_k\pi_{\theta_k}$. However, if $\mu$ denotes the Lebesgue measure on $\Theta$, and if $\Theta$ is a
given interval, say $\Theta=[a,b]$, then the distribution of the latent variable $t$ is spread out over this interval and
$\pi_f$ represents a \textit{continuous mixture}. In this paper we consider continuous mixtures and the problem of
identifying the mixing density $f$ when a sample of the continuous mixture
$\pi_f$ is observed. Note that when $t$ is
a location parameter, the problem of estimating $f$ is referred to as
a \emph{deconvolution problem}, which has received considerable attention in
the nonparametric statistics literature since
\cite{fan1991}. 

Continuous mixtures have been used in very numerous and various fields of
application.  We just give some recent examples to show that continuous
mixtures are still of much interest from an application point of view.  The
video-on-demand traffic can be modeled by a continuous Poisson mixture for the
purpose of efficient cache managing \citep{olmos14}. In time-resolved
fluorescence, where photon lifetimes have exponential distribution and
parameters depend on the emitting molecules, typically continuous mixtures of
exponential distributions are observed \citep{lakowicz99,valeur01}.  When $t$
is a scale parameter, the distribution $\pi_f$ is called a \emph{scale
  mixture}.  Scale mixtures of uniforms are also related to multiplicative
censoring introduced in \citet{vardi89} and length-biased data. A recent
application in nanoscience of the latter are length measurements of carbon
nanotubes, where observations are partly censored \citep{kvam08}.
Exponential mixtures play a significant role in natural sciences phenomena of
discharge or disexcitation as e.g. radioactive decays, the electric discharge
of a capacitor or the temperature difference between two objects.  Several
examples of applications of the exponential mixture model can be found in the
references of the seminal paper \cite{jewell82}.  

Not only for
applications, as well from a mathematical point of view, scale mixtures are
particularly interesting as they define classes of densities that verify some
monotonicity constraints. It is well known that any monotone non-increasing
density function with support in $(0,+\infty)$ can be written as a mixture of
uniform densities U$[0,t]$ \citep[][p. 158]{feller71}.  Moreover, a
\textit{$k$-monotone} density is defined as a non-increasing, convex density
function $h$ whose derivatives satisfy for all $j=1,\dots,k-2$ that
$(-1)^jh^{(j)}$ is non-negative, non-increasing and convex. One can show that
any $k$-monotone density can be represented by a scale mixture of Beta
distributions $B(1,k)$.  Furthermore, densities that are $k$-monotone for any
$k\geq1$, also called \textit{completely monotone} functions, can be written as
a continuous mixture of exponential distributions \citep{balabdaoui07}.

The literature provides various approaches for the estimation of the mixing density, as for
example the nonparametric maximum likelihood estimate (NPMLE). A characteristic
feature of this estimator is that it yields a discrete mixing distribution
\citep{laird78, lindsay83}. This appears to be unsatisfactory if we have
reasons to believe that the mixing density is indeed a smooth function. In this
case a functional approach is more appropriate, which relies on smoothness
assumptions on the mixing density $f$. In \citet{zhang90} kernel estimators are
constructed for mixing densities of a location parameter. \citet{goutis97}
proposes an iterative estimation procedure also based on kernel methods.  
\citet{asgharian12}  show strong uniform consistency of kernel estimators in the specific case of multiplicative censoring. 
In the same setting, \citet{andersen01} consider the linear operator $K$  verifying $\pi_f =Kf$ and estimate $f$ by an SVD reconstruction in the orthonormal basis of eigenfunctions of $K$. 
For
mixtures of discrete distributions, that is when $\pi_t$ are densities with
respect to a counting measure on a discrete space, orthogonal series estimators
have been developed and studied in \citet{hengartner97} and
\citet{roueff05}. For such mixtures, these estimators turn out to enjoy similar
or better rates of convergence than the kernel estimator presented in
\citet{zhang95}. 
\citet{comte14} present a projection estimator 
based on Laguerre functions that has the specific feature that the support of
the mixing density $f$ is not a compact as usual, but  the entire positive real
line. \citet{belomestny-schoemakers14}  extend the class of scale
mixtures and derive estimation methods based on the Mellin transform.

In this paper we show that orthogonal series estimators can be provided in a
very general fashion to estimate mixing densities with compact supports.  In
contrast to \citet{andersen01}, who consider only the case of scale mixtures of
uniforms, our approach applies to a large variety of continuous mixtures as our
numerous examples demonstrate.  In the {exponential mixture case}, in
particular, we exhibit an orthogonal series estimator achieving the minimax
rate of convergence in a collection of smoothness classes without requiring a
prior knowledge of the smoothness index. In other words, we provide an adaptive
estimator of the mixing density of an exponential mixture.

The paper is organized as follows. In Section \ref{sec_orthog_method} the
general construction of an orthogonal series estimator is presented and the
estimator is applied in several different mixture settings.  In Section
\ref{sec_properties} we derive upper bounds on the rate of convergence of the
mean integrated squared error of the estimator on some specific smoothness
classes.  In Section \ref{sec:appr-class} the approximation classes used for
the convergence rate are related to more meaningful smoothness classes defined
by weighted moduli of smoothness.  Section \ref{sec_lowerbound} is concerned
with the investigation of the minimax rate. On the one hand, a general lower
bound of the MISE is provided and on the other hand, some specific cases are
studied in detail.  Section \ref{sec_support} provides a consistent estimator
of the support of the mixing density. Finally, the performance of the
projection estimator is evaluated by a simulation study in different mixture
settings in Section \ref{sec_simul}. The Appendix provides some technical
results.

\section{Estimation Method}\label{sec_orthog_method}
In this section we develop an orthogonal series estimator and we provide
several examples, namely for mixtures of exponential, Gamma, Beta and uniform 
densities.

\subsection{Orthogonal Series Estimator}\label{sec_general_estim}
Throughout this paper the following assumption will be used.
\begin{ass}\label{ass_model}
  Let $\zeta$ be a dominating measure on the observation space
  $(\Xset,\Xfield)$.  Let $\{\pi_t,t\in\Theta\}$ be a
  parametric collection of densities with respect to $\zeta$. Furthermore, let
  the parameter space $\Theta=[a,b]$ be a compact interval with known endpoints
  $a<b$ in $\R$. We denote by $X, X_1,\dots,X_n$ an i.i.d. sample from the
  mixture distribution density $\pi_f$ defined by (\ref{def_mixture_dens}) with
  $\mu$ equal to the Lebesgue measure on $[a,b]$.
\end{ass}

For convenience, we also denote by $\pi_t$ and $\pi_f$ the probability measures
associated to these densities. Moreover we will use the functional analysis
notation $\pi_t(h)$ and $\pi_f(h)$, for the integral of $h$ with respect to
these probability measures.

The basic assumption of our estimation approach is that the mixing density $f$ in (\ref{def_mixture_dens}) is square
integrable, that is $f\in L^2[a,b]$. Then, for any complete orthonormal basis $(\psi_k)_{k\geq1}$ of the Hilbert space
$\H= L^2[a,b]$, the mixing density $f$ can be represented by the orthogonal series $f(t) = \sum_{k\geq1}
\coeff_k\psi_k(t)$, where the coefficients $\coeff_k$ correspond to the inner products of $f$ and $\psi_k$.  If we have
estimators $\hat \coeff_{n,k}$ of those coefficients, then an estimator of the mixing density $f$ is obtained by
$\sum_{k=1}^m\hat\coeff_{n,k}\psi_k.$

To construct estimators $\hat\coeff_{n,k}$, we remark that the following relation holds:  
Let $g$ be a nonnegative integrable function on $\R$. Define the function $\varphi$ on $[a,b]$ by the conditional expectations 
\begin{equation}\label{def_phi_g}
\varphi(t) = \pi_t(g) =
\int_{x\in\Xset} g(x)\pi_t(x)\;\zeta(\rmd x)\;,\quad t\in[a,b]\;.
\end{equation}
Suppose that $\varphi$ belongs to $\H$. The mean
$\pi_f(g)$ can be written as the inner product of $f$ and $\varphi$. Namely, by the definition of $\pi_f$ in
(\ref{def_mixture_dens}) and Fubini's theorem,
$$ 
\pi_f(g)=\int_{x\in\Xset} g(x) \pi_{f}(x) \;\zeta(\rmd x)
= \int_a^b f(t)\int_{x\in\Xset}g(x)\pi_t(x)\;\zeta(\rmd x)\,\rmd t 
= \langle f,\varphi\rangle_{\H}\;.
$$
Consequently, by the strong law of large numbers, $\frac1n\sum_ig(X_i)$ is a
consistent estimator of the inner product $\langle f,\varphi\rangle_{\H}$ based
on an i.i.d. sample $(X_1,\ldots,X_n)$ from the mixture density $\pi_{f}$
defined in (\ref{def_mixture_dens}).

We make the following assumption under which  the orthogonal series estimator makes sense.
\begin{ass}\label{ass_Vm}
  Assumption~\ref{ass_model} holds and there exists a sequence $(g_k)_{k\geq1}$
  of $\Xset\to\R$ functions such that $(\varphi_k)_{k\geq1}$ is a dense
  sequence of linearly independent functions in $\H$, where
  $\varphi_k(t)=\pi_{t}(g_k)$ as in~(\ref{def_phi_g}).
\end{ass}
We then proceed as follows.  Using linear combinations of the $\varphi_k$'s, a sequence of orthonormal functions
$\psi_1,\psi_2,\dots$ in $\H$ can be constructed, for instance by the Gram-Schmidt procedure. Say that $\psi_k$ writes as
$\sum_{j=1}^kQ_{k,j}\varphi_j$ with an array $(Q_{k,j})_{1\leq j\leq k}$ of real values that are computed beforehand. Then we
define estimators of $\coeff_k =\langle f, \psi_k\rangle_{\H} = \sum_{j=1}^kQ_{k,j}\langle f, \varphi_j\rangle_{\H} $ by the empirical means
\begin{equation*}
\hat\coeff_{n,k} = \frac 1n\sum_{i=1}^n\sum_{j=1}^kQ_{k,j}g_j(X_i)\;.
\end{equation*}
Finally, for any integer $m$,  an estimator of $f$ is given by 
\begin{equation}\label{def_proj_estim}
\hat f_{m,n} = \frac1n \sum_{k=1}^m \hat \coeff_{n,k}\psi_k =  \frac1n
\sum_{i=1}^n \sum_{j,k,l=1}^m Q_{k,j}Q_{k,l}g_j(X_i)\varphi_l\;,
\end{equation}
with the convention $Q_{k,j}=0$ for all $j>k$.
We refer to $\hat f_{m,n}$ as the \textit{orthogonal series estimator} or the \textit{projection estimator} of approximation order $m$. 

Define the subspaces 
\begin{equation}
  \label{eq:defVm}
  V_m =\text{span}(\varphi_1,\dots,\varphi_m)\;,\quad\text{ for all }m\geq1\;.
\end{equation} 
By Assumption~\ref{ass_Vm}, the sequence $(V_m)_m$ is strictly increasing,
$V_m$ has dimension $m$ for all $m$, and $\cup_mV_m$ has closure equal to
$\H$. By construction the orthogonal series estimator $\hat f_{m,n}$ belongs to
$V_m$. Consequently, the best squared error achievable by $\hat f_{m,n}$ is
$\|f-P_{V_m}(f)\|_{\H}^2$, where $\|\cdots\|_{\H}$ denotes the norm associated
to $\H$ and $P_{V_m}$ the orthogonal projection on the space $V_m$. Hence once
the functions $g_k$ are chosen, the definition of the subspaces $V_m$ follows
and the performance of the estimator will naturally depend on how well $f$ can
be approximated by functions in $V_m$.
% Note that it can be impossible to proceed conversely. That means
% that, if we first fix $V_m$, then it can occur that there are no adequate
% functions $g_k$ such that $\text{span}(\varphi_1,\dots,\varphi_m) = V_m$.
It is thus of interest to choose a sequence $(g_k)_{k\geq1}$ yielding a
meaningful sequence of approximation spaces $(V_m)_m$.  In the context of scale
family mixtures (but not only, see \citet{roueff05}), polynomial spaces appear
naturally.  Indeed, for any function $g$, we have $\pi_t(g)=\pi_1(g(t\cdot))$, so that,
provided that $\pi_1$ has finite moments, if $g$ is polynomial of degree $k$,
so is $\varphi(t)=\pi_t(g)$. The following assumption slightly extends this choice
for the two following reasons. First, a scale family is not always
parameterized by its scale parameter but by its inverse (as for the exponential
family).  Second, it will appear that the choice of $(g_k)_{k\geq1}$ not only
influences the approximation class (and thus the bias) but also the
variance. It may thus be convenient to allow the $g_k$'s not to be polynomials, while still remaining 
in the context of polynomial approximation.
This goal is achieved by the following assumption.

\begin{ass}\label{ass_polynom_L2ab}
  Assumption~\ref{ass_Vm} holds and there exist two real numbers $a'<b'$ and a
  linear isometry $T$ from $\H$ to $\H'=L^2[a',b']$ such that, for all
  $k\geq1$, $T\varphi_k$ is a polynomial of degree $k-1$. We denote by $T^{-1}$
  the inverse isometry.
\end{ass}

To compute the coefficients $Q_{k,j}$ under Assumption~\ref{ass_polynom_L2ab}, one may rely on the well known
Legendre polynomials which form an orthogonal sequence of
polynomials in $\H'=L^2[a',b']$. Indeed, by choosing
$g_k$ so that $T\varphi_k$ is the polynomial $t^{k-1}$, as will be illustrated in all the
examples below, the constants $Q_{k,j}$
are the coefficients of the normalized Legendre polynomials $\sum_{j=1}^kQ_{k,j}t^{j-1}$. 
Let us recall the definition of the Legendre polynomials.
\begin{definition}[Legendre polynomials]
  Let $a'<b'$ be two real numbers and denote $\mu=(a'+b')/2$ and
  $\delta=(b'-a')/2$. The \emph{Legendre polynomials} associated to the interval $[a',b']$
  are defined as the polynomials $r_k(t) = \sum_{l=1}^k
  R_{k,l} t^{l-1}$, where the coefficients $R_{k,l}$ are given by the following
  recurrence relation
\begin{align*}
R_{k+1,l} = R_{k,l-1} +\mu  R_{k,l}-\beta_k R_{k-1,l}\;,\quad \text{for all } k,l\geq1\;,
\end{align*}
with $ R_{1,1}=1$ and $R_{k,l}=0$ for all $l>k$, $\beta_1=2\delta$ and
$\beta_k=\delta^2(k-1)^2/(4(k-1)^2-1)$ for $k\geq2$. The obtained sequence
$(r_k)_{k\geq1}$ is orthogonal in $\H'=L^2([a',b'])$ with norms given by
$\|r_k\|_{\H'} = \sqrt{\beta_1\dots\beta_k}$. Hence, the coefficients of the
\emph{normalized Legendre polynomials} are defined by the relation
\begin{equation}\label{def_legendre_coeff}
Q_{k,l} = \frac{R_{k,l}}{\sqrt{\beta_1\dots\beta_k}}\;,\quad  \text{for all } k,l\geq1\;.
\end{equation} 
\end{definition}

\subsection{Examples}\label{subsec_orthog_examples}
 
For illustration we exhibit in this section the orthogonal series estimator in
some special cases. Some scale mixtures are presented. As an example for a
non scale mixture we also consider Gamma shape mixtures.

\noindent\textit{Example 1 (a). Exponential Mixture.}
We first consider  continuous exponential mixtures as they play a meaningful role in physics.  
That is, we consider $\pi_t(x) = t\rme^{-tx}$.  For the orthogonal series estimator we choose the functions
$g_k(x)=\1\left\{x>k-\frac12\right\}$ for $k\geq1.$ By (\ref{def_phi_g}), we
obtain 
\begin{equation}
  \label{eq:varphi_exp_a}
  \varphi_k(t) = \rme^{-(k-\frac12)t}\;.
\end{equation} 
We claim that the $\varphi_k$'s can be transformed into polynomials in the 
space $\H'=L^2[\rme^{-b},\rme^{-a}]$. Indeed, define, for all $f\in\H=L^2[a,b]$,
\begin{equation}\label{eq:T_exp_a}
Tf(t) = f(-\log t)/\sqrt t\;,\quad t\in[\rme^{-b},\rme^{-a}]\;.
\end{equation}
Then one has $\langle Tf, Tg\rangle_{\H'}=\langle f, g\rangle_{\H}$, hence $T$
is an  isometry  from $\H$ to $\H'$. Moreover  $T\varphi_k(t) = t^{k-1}$ are polynomials. Denote by $p_k(t) = \sum_{j=1}^kQ_{k,j}t^{j-1}$ the Legendre polynomials in $\H'$ with coefficients $Q_{k,j}$ defined by (\ref{def_legendre_coeff}) with $a'=\rme^{-b}$ and $b'=\rme^{-a}$. Denote by $T^{-1}$  the inverse operator of $T$  given by $T^{-1}h(t) = \rme^{-t/2}h(\rme^{-t})$. 
Since $T^{-1}$ is a linear isometry, we get that the functions $\psi_k \equiv T^{-1}p_k =  \sum_{j=1}^kQ_{k,j}\varphi_j$ are orthonormal in $\H$.
Consequently, an orthonormal series estimator is given by 
\begin{equation}\label{def_orthog_estim_3a}
\hat f_{m,n}(t) = \frac1n \sum_{k,j,l=1}^m\sum_{i=1}^n \1\left\{X_i>j-\frac12\right\} Q_{k,j}Q_{k,l}\rme^{-(l-\frac12)t}\;.
\end{equation}

\noindent\textit{Example 1 (b). Exponential Mixture.}
The choice of the functions $g_k$ is not unique and needs to be done with
care. For illustration, consider once again exponential mixtures with $\pi_t(x)
= t\rme^{-tx}$. 
%In Examples 1 and 2, we choose $g_k(x)=a_kx^{k-1}$, as we observe that in the context of a scale family, this leads to $\varphi_k(t)=t^{k-1}$ provided that $a_k=1/\int x^{k-1}\pi_1(\rmd x)$. Although the exponential family is parameterized by the inverse of the scale parameter, one can still choose $g_k$'s that are polynomials, leading to a completely different estimator from Example 3 (a).  
This time we take
$$
g_k(x) = a_k x^k \quad \text{ with } a_k = \left(\int x^{k}\pi_1(\rmd x)\right)^{-1}=1/k!$$
and hence $\varphi_k(t) = t^{-k},$ for $k\geq1$. 
To relate $\varphi_k$ to polynomials, define the isometry $\tilde T$ from $\H$
to $\tilde\H = L^2[1/b,1/a]$ by $\tilde Tf(t) = \frac1tf(\frac 1t)$. We
have $\tilde T\varphi_k(t) =t^{k-1}$ for all $k\geq1$. Furthermore, denote by
$\tilde T^{-1}$ the inverse of $\tilde T$ satisfying $\tilde
T^{-1}h=\frac1th(\frac1t)$.  Let $\tilde p_k(t) = \sum_{j=1}^k Q_{k,j}t^{j-1}$
be the Legendre polynomials in $\tilde\H$ defined with $a'=1/b$ and $b'=1/a$. Since $\tilde T^{-1}$
is an isometry, $\psi_k \equiv \tilde T^{-1}p_k = \sum_{j=1}^k 
Q_{k,j}\varphi_{j}$ are orthonormal functions in $\H$ and the orthonormal
series estimator is given by
$$
\hat f_{m,n}(t) = \frac1n \sum_{k,j,l=1}^m\sum_{i=1}^n \frac{  Q_{k,j}Q_{k,l}X_i^j}{j!}t^{-l}\;.
$$

\noindent\textit{Example 2. Gamma Shape Mixture.}
Polynomial estimators can be used in the context where the mixed parameter is
not necessarily a scale parameter. As pointed out earlier, they have first been
used for mixtures on a discrete state space $\Xset$, such as Poisson mixtures,
see \cite{hengartner97} and \cite{roueff05}. Let us consider the Gamma shape mixture
model.  Parametric Gamma shape mixtures have been considered in
\cite{venturini:2008}. For this model $\pi_t$ is the Gamma density with shape
parameter $t$ and a fixed scale parameter (here set to 1 for simplicity),
$$
\pi_t(x)=\frac{x^{t-1}}{\Gamma(t)}\;\rme^{-x},\quad t\in[a,b]\;,
$$
where $\Gamma$ denotes the Gamma function. This model has a continuous state
space ($\zeta$ is the Lebesgue measure on $\R_+$) and is not a scale 
mixture. We shall construct $g_k$ and $\varphi_k=\pi_{\cdot}(g_k)$ such that
Assumption~\ref{ass_polynom_L2ab} holds with $T$ being the identity and
$\varphi_k(t)=t^{k-1}$.  Consider the following sequence of polynomials,
$p_1(t)=1$, $p_2(t)=t$, ..., $p_k(t)=t(t+1)\dots(t+k-2)$ for all $k\geq2$.
Since $(p_k)_{k\geq1}$ is a sequence of polynomials with degrees $k-1$, there
are  coefficients $(\tilde{c}_{k,l})_{1\leq l\leq k}$ such that
$t^{k-1}=\sum_l\tilde{c}_{k,l}p_l(t)$ for $k=1,2,\dots$.  A simple recursive
formula for computing $(\tilde{c}_{k,l})_{1\leq l\leq k}$ is provided in
Lemma~\ref{lem:rec_gamma_pol} in the Appendix, see
Eq.~(\ref{eq:rec_gamma_pol}).  Observe that, for any $l\geq1$,
$$
\int x^{l-1}\;\pi_t(x)\rmd x =\frac{\Gamma(t+l-1)}{\Gamma(t)}=p_l(t)\;.
$$
Hence, setting $g_k(x)=\sum_l\tilde{c}_{k,l}x^{l-1}$, we obtain
$$
\varphi_k(t)=\pi_t(g_k)=\sum_l\tilde{c}_{k,l}p_l(t)=t^{k-1}\;,
$$
and thus Assumption~\ref{ass_polynom_L2ab} holds with $T$ being the identity
operator and $\varphi_k(t)=t^{k-1}$. Define $(Q_{k,l})_{k,l}$ as the coefficients
of Legendre polynomials on $\H=L^2([a,b])$, that is as
in~(\ref{def_legendre_coeff}) with $a'=a$ and $b'=b$. The 
polynomial estimator defined by~(\ref{def_orthog_estim_3a}) reads
$$
\hat f_{m,n}(t) = \frac1n \sum_{i=1}^n\sum_{k,j,l=1}^m Q_{k,j} Q_{k,l}
\sum_{h=1}^j \tilde{c}_{j,h}X_i^{h-1}t^{l-1}\;.
$$

\noindent\textit{Example 3. Scale Mixture of Beta Distributions or Uniform Distributions.}
It is well known that any $k$-monotone density, for $k\geq1$, can be represented by a scale mixture of Beta distributions $B(1,k)$  \citep{balabdaoui07} with 
$$\pi_t(x) = \frac kt\left(1-\frac xt\right)^{k-1},\quad\text{ for } x\in[0,t]\;.$$
Note that if $k=1$, then $\pi_t$ is the uniform density $U(0,t)$. 
We take
$$
g_p(x) = a_p x^{p-1} \quad \text{ with } a_p = \left(\int x^{p-1}\pi_1(\rmd x)\right)^{-1} =\frac1{k~\beta(p,k)},\quad p\geq1\;,
$$
where $\beta(a,b) = \int_0^1 t^{a-1}(1-t)^{b-1}\rmd t$ denotes the Beta
function.  It follows that $\varphi_p(t)=t^{p-1}$.  As in the preceding
example, if $f\in\H$ then an orthogonal series
estimator $\hat f_{m,n}$ of $f$ can be constructed by using Legendre 
polynomials $p_k(t) = \sum_{j=1}^kQ_{k,j}t^{j-1}$ where the coefficients $Q_{k,j}$ are defined as in~(\ref{def_legendre_coeff})  with $a'=a$ and $b'=b$. Then according to
(\ref{def_proj_estim}), the corresponding orthogonal series estimator is given
by
$$
\hat f_{m,n}(t) = \frac1n \sum_{j,p,l=1}^m
\sum_{i=1}^n
Q_{p,j}Q_{p,l}\frac{X_i^{j-1}}{k\beta(j,k)}t^{l-1}\;.
$$
In Example 1 (b) we considered the same functions $g_p$ but here
Assumption~\ref{ass_polynom_L2ab} holds with $T$ equal to the identity operator
on $\H$.  This difference relies on the parametrization of the exponential
family by the inverse of the scale parameter.
%In Examples 1 and 2, we choose $g_k(x)=a_kx^{k-1}$, as we observe that in the context of a scale family, this leads to $\varphi_k(t)=t^{k-1}$ provided that $a_k=1/\int x^{k-1}\pi_1(\rmd x)$. Although the exponential family is parameterized by the inverse of the scale parameter, one can still choose $g_k$'s that are polynomials, leading to a completely different estimator from Example 3 (a).  

\noindent\textit{Example 4. Mixture of exponential distributions with location parameter.} 
The estimator also applies to the deconvolution setting. As an example,
consider $X=Y+\theta$ where $Y$ and $\theta$ are independent random variables,
$Y$ has exponential distribution with mean 1 and $\theta$ has unknown density
$f$ supported on $[a,b]$.  The density of $X$ is given by $\pi_f (x) =\int_a^b
\pi_t(x)f(t)\rmd t$ with $\pi_t(x) = \rme^{-(x-t)}\1\{x>t\}$. Let $g_1(x)=1$
and 
\begin{equation}
  \label{eq:gk-location-exp}
g_k(x)=x^{k-1}-(k-1)x^{k-2}\;,  
\end{equation}
for $k\geq2$. Then $\varphi_k(t) = t^{k-1}$
for $k\geq1$.  The estimator $\hat f_{m,n}$ of $f$ is then given by
\begin{align*}\label{estim_exp_deconv}
\hat f_{m,n}(t) = \frac1n\sum_{i=1}^n\sum_{j,k,l=1}^mQ_{k,l}Q_{k,j}g_j(X_i)t^{l-1}\;,
\end{align*}
where $Q_{k,l}$ are the Legendre coefficients defined by~(\ref{def_legendre_coeff})   with $a'=a$ and $b'=b$.

\section{Analysis of the Orthogonal Series Estimator}\label{sec_properties}

In this section the properties of the orthogonal series estimator are
analyzed. % First, we compute the bias and the variance of $\hat f_{m,n}$. Then
% upper and lower bounds of the mean integrated squared error (MISE) are
% developed. Finally, we will show that in the case of an exponential mixture,
% the orthogonal series estimator defined by (\ref{def_orthog_estim_3a}) achieves
% the minimax rate on some specific smoothness classes of functions.

\subsection{Bias, Variance and MISE}
It is useful to write the orthogonal series estimator $\hat f_{m,n}$ defined in (\ref{def_proj_estim}) in matrix notation. Therefore, we introduce the $m\times m$--matrix $Q=(Q_{k,j})_{k,j}$, where $Q_{k,j}=0$ for all $j>k$, and the $m$--vectors 
\begin{align*}
&\Phi = [\varphi_1,\ldots,\varphi_m]^T\;,\quad 
\Psi = [\psi_1,\dots,\psi_m]^T = Q\Phi\;,\\ 
&{\bf g}(x) = [g_1(x),\ldots,g_m(x)]^T\;,\quad
{\bf\hat g} = \frac1n\sum_{i=1}^n{\bf g}(X_i)\;,\\
&{\bf \coeff} = [\coeff_1,\dots,\coeff_m]^T = \langle\Psi,f\rangle_{\H}\;,\quad 
{\bf \hat \coeff} = [\hat\coeff_{n,1},\dots,\hat\coeff_{n,m}]^T = Q{\bf \hat g}\;.
\end{align*}
It follows that the orthogonal series estimator can be written as
$$
\hat f_{m,n} = {\bf \hat \coeff}^T\Psi ={\bf\hat g}^TQ^TQ\Phi\;.$$ 
Further, let $\Sigma = \pi_f({\bf g}{\bf g}^T)-\pi_f({\bf g})\pi_f({\bf
  g})^T$ be the covariance matrix of ${\bf g}(X_1)$. The MISE is defined by
$\E\left\|\hat f_{m,n}-f\right\|_{\H}^2$.   
The orthogonal projection of $f$ on $V_m$ is denoted by 
$$P_{V_m}f = {\bf \coeff}^T\Psi=\sum_{k=1}^m\coeff_{n,k}\psi_k\;.$$ 

It is clear that the orthogonal series estimator $\hat f_{m,n}$ is an unbiased
estimator of $P_{V_m}f$. Furthermore, by the usual argument, the MISE is
decomposed into two terms representing the integrated variance and integrated
squared bias, as summarized in the following result, whose proof is standard
and thus omitted. 

\begin{prop}\label{prop_MISE}
  Suppose that Assumption~\ref{ass_Vm} holds.  The orthogonal series estimator
  $\hat f_{m,n}$ defined in (\ref{def_proj_estim}) satisfies
\begin{enumerate}[(i)]
	\item For every $t\in[a,b]$, $~~\E[\hat f_{m,n}(t)] = P_{V_m}f(t).$
	\item \label{prop_item_var} For every $t\in[a,b]$, $~~\mathrm{Var}(\hat f_{m,n}(t)) = \frac1n\Psi^T(t) Q \Sigma Q^T\Psi(t).$
	\item \label{prop_item_mise} $\E\left\|\hat f_{m,n}-f\right\|_{\H}^2$ =
          $\|P_{V_m}f-f\|_{\H}^2 + \frac1n \mathrm{tr}\left(Q\Sigma Q^T\right).$
\end{enumerate}
\end{prop}

% \begin{proof}
% \begin{enumerate}[(i)]
% 	\item Clear, since $\E[{\bf \hat\coeff}]={\bf \coeff}$. 
% 	\item For the variance we find 
% 	$$\mathrm{Var}(\hat f_{m,n}(t)) 
% 	= \E[\hat f_{m,n}(t)^2] - \E[\hat f_{m,n}(t)]^2 \\
% 	= \Psi(t)^T Q \E[{\bf\hat g}{\bf\hat g}^T]Q^T\Psi(t) - (P_{V_m}f)^2\;. $$
% 	We see that 
% 	\begin{align*}
% 	\E[{\bf\hat g}{\bf\hat g}^T] 
% 	&= \frac1n\E\left[{\bf g}(X){\bf g}(X)^T\right] +\frac{n-1}n \E [{\bf g}(X)]\E [{\bf g}(X)^T]
% 	= \frac1n\Sigma + \langle \Phi,f\rangle\langle \Phi,f\rangle^T
% 	\end{align*}
% 	and $\Psi^TQ\langle \Phi,f\rangle  = \Psi^T  \langle\Psi,f\rangle={\bf \coeff}^T\Psi =P_{V_m}f $. This implies $(ii)$.
% 	\item By the Pythagorean theorem
% $$\|\hat f_{m,n}-f\|_{\H}^2 = \|\hat f_{m,n}-P_{V_m}f\|_{\H}^2+\|P_{V_m}f-f\|_{\H}^2\;.$$
% Moreover, by Fubini's thoerem
% \begin{align*}
% \E\|\hat f_{m,n}-P_{V_m}f\|_{\H}^2
% &= \int_a^b  \mathrm{Var}(\hat f_{m,n}(t))\rmd t 
% = \frac 1n\int_a^b  \mathrm{Var}(\hat f_{m,1}(t))\rmd t \;.
% \end{align*}
% \end{enumerate}
% This concludes the proof.
% \end{proof}

An important issue for orthogonal series estimators $\hat f_{m,n}$ is the choice of the approximation order $m$. % The
% stability of solving $\pi_f(x) = \int \pi_t(x)f(t)\rmd t$ is controlled by choosing $m$ to be small relative to $n$.  When
% $m$ is chosen rather large, then the number of coefficients $\coeff_k$ to be estimated is large, which entails a large
% variance of the estimator. However, when $m$ is too small, the bias gets large, since $\hat f_{m,n}$ is an estimator of the
% orthogonal projection of $f$ on $V_m$ and $V_m\subset V_{m+1}$. Hence, the problem is to choose $m$ such that a bias-variance
% trade-off is achieved.
The integrated squared bias $\|P_{V_m}f-f\|_{\H}^2$ only depends on how
well $P_{V_m}f$ approximates $f$, whose rate of convergence depends on the smoothness class to which belongs the
density $f$. To be more precise, define for any
approximation rate index $\alpha$ and radius $C$, the approximation class
\begin{equation}\label{def_fnspace_Cuc}
\mathcal{C}(\alpha,C) = \{f\in\H : \|f\|_{\H}\leq C \text{ and } \|P_{V_m}f-f\|_{\H}\leq C\,m^{-\alpha} \text{ for all } m\geq 1\}\;.
\end{equation}
So when the mixing density $f$ belongs to $\mathcal{C}(\alpha,C)$, then the
bias of the orthogonal series estimator $\hat f_{m,n}$ is well controlled,
namely it decreases at the rate $m^{-\alpha}$ as $m$ increases.  Furthermore,
denote the set of densities in $\H$ by $\H_1 = \{f\in\H: f\geq0, \int_a^b
f(t)\rmd t = 1\}.$ We will investigate the rate of convergence of $\hat
f_{m,n}$ in $\H$ when $f\in\mathcal{C}(\alpha,C)\cap\H_1$.  We will obtain
the best achievable rate in the case of exponential mixtures and almost
the best one in the case of Gamma shape mixtures.

\subsection{Upper Bound of the MISE}\label{subsec_upperbound}
We now provide an upper bound of the MISE for the orthogonal series estimator
based on Legendre polynomials, that is, when Assumption~\ref{ass_polynom_L2ab}
holds.

To show an upper bound of the MISE we use the following property
\citep[see][Lemma A.1]{roueff05}.  If $\lambda > \frac{2+a'+b'}{b'-a'} +
\sqrt{1+\frac{2+a'+b'}{b'-a'}}$, then the coefficients of the normalized
Legendre polynomials in $L^2[a',b']$ defined by (\ref{def_legendre_coeff})
verify
\begin{equation}\label{rel_legendrecoeff}
\sum_{l=1}^kQ_{k,l}^2 = O(\lambda^{2k})\;,\quad\text{as } k\to\infty\;.
\end{equation}

% \begin{lem}\label{lem_upperbound_var} 
% Let the coefficients $Q_{k,j}$ of the orthogonal series estimator $\hat f_{m,n}$ defined in (\ref{def_proj_estim})  be  the coefficients of the normalized Legendre polynomials in some space $L^2[a',b']$ given by (\ref{def_legendre_coeff}). Let $\lambda > \frac{2+a'+b'}{b'-a'} + \sqrt{1+\frac{2+a'+b'}{b'-a'}}$. 
% If the variances $\mathrm{Var}(g_k(X))$ are bounded by a common constant $C_0$ for all $k\geq1$, then 
% $$ \int_a^b  \mathrm{Var}(\hat f_{m,1}(t))\rmd t=O\left(m\lambda^{2m}\right)\;,\quad \text{ as }m\to\infty\;.$$ 
% \end{lem}

% \begin{proof}
% Note that from the condition that all variances $\mathrm{Var}(g_k(X))$ are bounded by a common constant $C_0$ it follows by the Cauchy-Schwarz inequality that all covariances $\Sigma_{k,l}=\mathrm{Cov}(g_k(X),g_l(X))$ are bounded by $C_0$.
% By Proposition \ref{prop_MISE} and by linearity of the trace, we obtain
% \begin{align*}
% \int_a^b  \mathrm{Var}(\hat f_{m,1}(t))\rmd t
% &= \int_a^b \mathrm{tr}\left(\Psi^T(t)Q\Sigma Q^T\Psi(t)\right)\rmd t
% = \mathrm{tr}\left(Q\Sigma Q^T\int_a^b \Psi(t)\Psi^T(t)\rmd t\right)\\
% &= \mathrm{tr}\left(Q\Sigma Q^T\right)
% = \sum_{k=1}^m\sum_{j=1}^m\sum_{l=1}^m Q_{k,j}Q_{k,l}\Sigma_{j,l} \\
% &\leq C_0 \sum_{k=1}^m \left( \sum_{j=1}^m |Q_{k,j}|\right)^2
% \leq mC_0 \sum_{k=1}^m  \sum_{j=1}^m |Q_{k,j}|^2\\
% &\leq K'm \sum_{k=1}^m\lambda^{2k}
% \leq Km\lambda^{2m}\;,
% \end{align*}
% where the first inequality comes from  (\ref{rel_legendrecoeff}) and $K$ and $K'$ are positive numbers.  This concludes the proof.
% \end{proof}
By combining Proposition \ref{prop_MISE} (\ref{prop_item_mise}) and 
the bound given in~(\ref{rel_legendrecoeff}) along with a normalization condition on the
$g_k$'s (Condition~(\ref{cond_bounded_var}) or Condition~(\ref{cond_bounded_var_caseb}) below), %Lemma \ref{lem_upperbound_var} 
we obtain the following asymptotic upper bounds of the MISE.
\begin{theo}\label{theo_mise_upperbound}
  Let $\alpha$ be a positive rate index and $C$ be a positive radius. Suppose
  that Assumption~\ref{ass_polynom_L2ab}
  holds with $f\in \mathcal{C}(\alpha,C)\cap \H_1$. Let $\hat f_{m,n}$ be defined
  by~(\ref{def_proj_estim}) with Legendre polynomials coefficients $Q_{k,j}$
  given by~(\ref{def_legendre_coeff}). Then the two following assertions hold.
  \begin{enumerate}[(a)]
  \item\label{item:case_a} If, for some  constants $C_0>0$ and $B\geq1$, we have
\begin{equation}\label{cond_bounded_var}
\mathrm{Var}(g_k(X))<C_0\,B^{2k}\quad\text{ for all } k\geq1\;.
\end{equation}  
Set $m_n=A\log n$  with 
\begin{equation}
  \label{eq:Aub}
A<\frac12\left\{\log B+\log\left(\frac{2+a'+b'}{b'-a'} + \sqrt{1+\frac{2+a'+b'}{b'-a'}}\right)\right\}^{-1}\;.  
\end{equation}
Then, as $n\to\infty$,
\begin{equation}\label{upperbound}
 \E\left\|\hat f_{m_n,n}-f\right\|_{\H}^2\leq C^2 m_n^{-2\alpha}(1+o(1))\;,
\end{equation}
where the $o$-term only depends on the constants $\alpha$, $C$, $a'$, $b'$, $A$
and $C_0$.
\item\label{item:case_b}  If, for some  constants $C_0>0$ and $\eta>0$, we have
\begin{equation}\label{cond_bounded_var_caseb}
\mathrm{Var}(g_k(X))<C_0\;k^{\eta k}\quad\text{ for all } k\geq1\;.
\end{equation}  
Set $m_n=A\log n/\log\log n$  with $A<\eta^{-1}$. 
Then, as $n\to\infty$,
\begin{equation}\label{upperbound_caseb}
 \E\left\|\hat f_{m_n,n}-f\right\|_{\H}^2\leq C^2 m_n^{-2\alpha}(1+o(1))\;,
\end{equation}
where the $o$-term only depends on the constants $\alpha$, $C$, $a'$, $b'$, $A$
and $C_0$.
  \end{enumerate}
\end{theo}
\begin{remark}
  \label{rem:variance_B}
  The larger $A$, the lower the upper bound in~(\ref{upperbound}). Hence, since
  $a',b'$ and $B$ directly depend on the $g_k$'s, the constraint~(\ref{eq:Aub})
  on $A$ indicates how appropriate the choice of the $g_k$'s is.
\end{remark}
\begin{remark}
  \label{rem:minimax-rate}
  In the examples treated in this paper, $C_0$ and $B$ or $\eta$ can be chosen
  independently of $f\in\mathcal{C}(\alpha,C)\cap \H_1$. Consequently, the
  bounds given in~(\ref{upperbound}) and~(\ref{cond_bounded_var_caseb}) show
  that $\hat f_{m_n,n}$ achieves the MISE rates $(\log n)^{-2\alpha}$ and
  $(\log(n)/\log\log n)^{-2\alpha}$, respectively, uniformly on
  $f\in\mathcal{C}(\alpha,C)\cap\H_1$. In the exponential mixture case, we show
  below that $\hat f_{m_n,n}$ of Example~1(a) is minimax rate adaptive in these
  classes (since $m_n$ does not depend on $\alpha$). In
  the Gamma shape mixture case, we could only show that $\hat f_{m_n,n}$ of
  Example~2 is minimax rate adaptive in these classes up to the multiplicative
  $\log\log n$ term.
\end{remark}
\begin{proof}
We first consider  Case~(\ref{item:case_a}).
By ~(\ref{eq:Aub}), we may choose a number $\lambda$ strictly lying between 
$\frac{2+a'+b'}{b'-a'} +
\sqrt{1+\frac{2+a'+b'}{b'-a'}}$ and $\rme^{1/(2A)}/B$.
Note that from  Condition~(\ref{cond_bounded_var}), it follows by the
Cauchy-Schwarz inequality that $|\Sigma_{k,l}|=|\mathrm{Cov}(g_k(X),g_l(X))|\leq
C_0B^kB^l$ for all $k,l$. Thus, we obtain
\begin{align*}
%\int_a^b  \mathrm{Var}(\hat f_{m,1}(t))\rmd t
% &= \int_a^b \mathrm{tr}\left(\Psi^T(t)Q\Sigma Q^T\Psi(t)\right)\rmd t
% = \mathrm{tr}\left(Q\Sigma Q^T\int_a^b \Psi(t)\Psi^T(t)\rmd t\right)\\
\mathrm{tr}\left(Q\Sigma Q^T\right)
&\leq C_0\sum_{k=1}^m\sum_{j=1}^k\sum_{l=1}^k |Q_{k,j}Q_{k,l}|\;B^jB^l \\
&\leq C_0 \sum_{k=1}^m \left( \sum_{j=1}^k Q_{k,j}^2\sum_{j=1}^k B^{2j}\right)\\
%\leq mC_0 \sum_{k=1}^m  \sum_{j=1}^m |Q_{k,j}|^2\\
%&\leq K'm \sum_{k=1}^m\lambda^{2k}
&\leq Km\{B\lambda\}^{2m}\;,
\end{align*}
where the last inequality comes from  (\ref{rel_legendrecoeff}) and $K$ is a
positive constant (the multiplicative term $m$ is necessary only for $B=1$).  
It follows by the decomposition of the MISE in Proposition \ref{prop_MISE} (\ref{prop_item_mise}) that
\begin{align*}
\E\left\|\hat f_{m_n,n}-f\right\|_{\H}^2
&\leq C^2m_n^{-2\alpha} +  Kn^{-1}m_n(B\lambda)^{2m_n}\\
&\leq C^2m_n^{-2\alpha}\left(1+ \frac{K}{C^2}n^{-1}m_n^{2\alpha+1}(B\lambda)^{2m_n} \right)\;.
\end{align*}
Now we have for $m_n = A\log n$ that
$$n^{-1}m_n^{2\alpha+1}(B\lambda)^{2m_n} = A^{2\alpha+1}(\log n)^{2\alpha+1} n^{2A\log B\lambda-1} = o(1)\;,$$
since $A<1/(2\log B\lambda)$.

Let us now consider Case~(\ref{item:case_b}). Proceeding as above, for any
$\lambda > \frac{2+a'+b'}{b'-a'} +\sqrt{1+\frac{2+a'+b'}{b'-a'}}$, we get
$\mathrm{tr}\left(Q\Sigma Q^T\right)\leq K\,C_0\lambda^{2m}m^{1+\eta m}$,
which yields
$$
\E\left\|\hat f_{m_n,n}-f\right\|_{\H}^2 \leq C^2m_n^{-2\alpha}\left(1+
  \frac{K}{C^2}n^{-1}m_n^{2\alpha+1+\eta m_n}\lambda^{2m_n} \right)\;.
$$
To conclude, it suffices to check that the log of the second term between
parentheses tends to $-\infty$ as $n\to\infty$ for $m_n=A\log n/\log\log n$
with $A<\eta^{-1}$, which is easily done.
\end{proof}

Let us check the validity of Condition (\ref{cond_bounded_var}) or
Condition~(\ref{cond_bounded_var_caseb}) for the above examples.

\noindent\textit{Example 1 (a). Exponential Mixture (continued).}
Condition (\ref{cond_bounded_var}) immediately holds with $B=C_0=1$ for the
exponential mixture of Example 1(a) since
$g_k(x)=\1\left\{x>k-\frac12\right\}$.

\noindent\textit{Example 1 (b). Exponential Mixture (continued).}
Interestingly, Condition (\ref{cond_bounded_var}) does not hold for Example~1(b), where a different choice of $g_k$'s is proposed. In fact, one finds that
$\log\mathrm{Var}(g_k(X))$ is of order $k\log(k)$. Hence, only
Condition~(\ref{cond_bounded_var_caseb}) holds and we fall in 
case~(\ref{item:case_b}) of Theorem~\ref{theo_mise_upperbound}.  Since a slower
rate is achieved in this case, this clearly advocates to choose the estimator
obtained in Example 1(a) rather than the one in Example 1(b) for the
exponential mixture model.

\noindent\textit{Example 2. Gamma Shape Mixture (continued).} We recall that here we set
$g_k(t)=\sum_{l=1}^k \tilde{c}_{k,l}t^{l-1}$ where the coefficients
$(\tilde{c}_{k,l})$ are those defined and computed in
Lemma~\ref{lem:rec_gamma_pol} of the Appendix. Using the
bound given by~(\ref{eq:rec_gamma_pol_sup_bound}) in the same lemma, we obtain that
$g_k(x)\leq k!(1\vee|x|^{k-1})$. It follows that $\pi_t(g_k^2)\leq (k!)^2
(1+\Gamma(t+2k-2)/\Gamma(t))$, and, for any $f\in\H_1$, $\pi_f(g_k^2)\leq
(k!)^2(1+\Gamma(b+2k-2)/\Gamma(b))$. Hence, by Stirling's formula, we find that
Condition~(\ref{cond_bounded_var_caseb}) holds for $\eta=4$ and some $C_0$
independent of $f\in\H_1$. % In fact the obtained rate will be shown to be
% minimax optimal in this case (see Theorem~\ref{theo_minimax_gammamix}). 

\noindent\textit{Example 3. Scale Mixture of Beta Distributions or Uniform Distributions (continued).}
We now verify Condition (\ref{cond_bounded_var})  for Beta mixtures and the
$g_p$ of Example 3. Note that we can write $X=\theta X_0$ with independent random variables $\theta\sim f$ and $X_0\sim
B(1,k)$. We have for all $p\geq1$
\begin{align*}
\mathrm{Var}(g_p(X))
&\leq %\E[g_p^2(X)] 
 \frac{\E[X^{2p-2}]}{k^2\beta^2(p,k)} 
= \frac{\E[\theta^{2p-2}]\E[X_0^{2p-2}]}{k^2\beta^2(p,k)}
\leq \frac{b^{2p-2} }{k^2\beta^2(p,k)}
\leq b^{2p-2} k^{2p-2} \;.
%\sim 2^{-k}(k!)^{-1} b^{2p-2} p^{k}\;,
%=\frac1{k(\beta(p,k))^2} \int_a^bf(t)\int_0^t\frac{x^{2p-2}}t\left(1-\frac xt\right)^{k-1}\rmd x \rmd t\\
%&= \frac{\beta(2p-1,k)}{k(\beta(p,k))^2} \int_a^bf(t)t^{2p-2}\rmd t
%\leq \frac{b^{2p-2}}{(\Gamma(p))^2\Gamma(k)}\;. FAUTE!!!
\end{align*}
%where the asymptotic equivalence holds as $p\to\infty$.
Hence Condition~(\ref{cond_bounded_var}) holds with $B=k$ if $b<1$, with   $B=bk$ if $b\geq1$. 

%\noindent\textit{Example 1. Uniform Mixture (continued).}
%In the setting of uniform mixtures where $g_k(x) = kx^{k-1}$ (Example
%1), where $X=\theta X_0$ with $\theta\sim f$ and $X_0\sim U[0,1]$, we have
%\begin{align*}
%\mathrm{Var}(g_k(X))
%\leq \E[k^2X^{2k-2}] = k^2\E[\theta^{2k-2}] \E[X_0^{2k-2}]
%\leq \frac{k^2b^{2k-2}}{2k-1}\;.
%\end{align*}
%Then (\ref{cond_bounded_var}) holds with $B$ given as in Example 3.

A close inspection of Example 3 indicates that it is a particular case
of the following more general result concerning mixtures of compactly supported
scale families.

\begin{lem}\label{lem:comp-support}
Suppose that Assumption~\ref{ass_model} holds in the context of a scale
mixture  on $\R_+$, that is, $\zeta$ is the Lebesgue measure on $\R_+$ and
$\pi_t=t^{-1}\pi_1(t^{-1}\cdot)$ for all $t\in\Theta=[a,b]\subset(0,\infty)$. 
Assume in addition that $\pi_1$ is compactly supported in $\R_+$.
Define, for all $k\geq1$, 
$$
g_k(x)=\left(\int x^{k-1}\pi_1(x)\,\rmd x\right)^{-1}\;x^{k-1}\;. 
$$
Then Assumption~\ref{ass_Vm} holds with $\varphi_k(t)=t^{k-1}$, and thus also does
Assumption~\ref{ass_polynom_L2ab} with $T$ being the identity operator on $L^2([a,b])$.
Moreover there exists $C_0$
and $B$ only depending on $\pi_1$ and
$b$ such that Condition~(\ref{cond_bounded_var}) holds.
\end{lem} 
\begin{proof}
  Using the assumptions on $\pi_1$ and Jensen's inequality, we have
\begin{equation*}%\label{cond_exp_mement}
B_1^m\leq \int x^m\;\pi_1(x)\,\rmd x \leq B_2^m\quad\text{for all
  $m\geq1$}\;,  
\end{equation*}
with $B_1=\int x\;\pi_1(x)\,\rmd x$ and $B_2>0$ such that the support of $\pi_1$ is included in $[0,B_2]$. 
The result then follows from the same computations as in Examples 3.
\end{proof}

An immediate consequence of Theorem~\ref{theo_mise_upperbound}
and Lemma~\ref{lem:comp-support} is the following.
\begin{cor}\label{cor:comp-support}
Under the assumptions of Lemma~\ref{lem:comp-support}, the estimator
$\hat f_{m,n}$ defined by~(\ref{def_proj_estim}) with Legendre polynomials coefficients
$Q_{k,j}$ given by~(\ref{def_legendre_coeff}) achieves the MISE rate
$(\log n)^{-2\alpha}$ uniformly on $f\in\mathcal{C}(\alpha,C)\cap\H_1$ for any
$\alpha>0$ and $C>0$.
\end{cor}

\textit{Example 4. Exponential mixture with location parameter (continued).} 
One can show that Condition (\ref{cond_bounded_var_caseb}) of Theorem
\ref{theo_mise_upperbound} is satisfied, so that the rate of the MISE of the estimator is $(\log (n)/\log\log n)^{-2\alpha}$.
 Indeed,
$$\E[X^r]%=\int x^r\pi_f(x)\rmd x 
= r!\int_a^bf(t)\sum_{j=0}^r\frac{t^j}{j!}\rmd t\leq r!\sum_{j=0}^r\frac{b^j}{j!}\leq r!\rme^b\;,$$
and thus, using the definition of $g_k$ in~(\ref{eq:gk-location-exp}), $\mathrm{Var}(g_k(X))
\leq 2(2k-2)!\rme^b 
\approx 2\sqrt{2\pi}\rme^{b-2k+2}(2k-2)^{2k-3/2}$.

\section{Approximation Classes}\label{sec:appr-class}
Although the approximation classes $\mathcal{C}(\alpha,C)$ appear naturally when studying the
bias  of the orthogonal series estimator defined in~(\ref{def_proj_estim}), it is
legitimate to ask whether such classes can be interpreted in a more intuitive
way, say using a smoothness criterion. This section provides a positive answer
to this question. 

\subsection{Weighted Moduli of Smoothness}
Let us recall the concept of weighted moduli of smoothness as introduced by
\citet{ditzian87} for studying the rate of polynomial approximations. For $a<b$
in $\R$, $f:[a,b]\to\R$, $r\in \N^*$ and $h\in\R$ denote by
$\Delta_{h}^r(f,\cdot)$ the \textit{symmetric difference }of $f$ of order $r$
with step $h$, that is
\begin{equation}\label{eq:diff_order_r}
\Delta_{h}^r(f,x)= \sum_{i=0}^r\left(r\atop i\right)(-1)^i f(x+(i-r/2)h)\;.
\end{equation}
with the convention that $\Delta_h^r(f,x)=0$ if $x\pm mh/2\notin[a,b]$.  Define
the step-weight function $\varphi$ on the bounded interval $[a,b]$ as
$\varphi(x)=\sqrt{(x-a)(b-x)}$.  Then for $f:[a,b]\to\R$ the \textit{weighted
  modulus of smoothness }of $f$ of order $r$ and with the step-weight function
$\varphi$ in the $L^p([a,b])$ norm is defined as
$$\omega_{\varphi}^r(f,t)_{p} = \sup_{0<h\leq t}\|\Delta_{h\varphi(\cdot)}^r(f,\cdot)\|_p\;.$$

We recall an equivalence relation of the modulus of smoothness with the so-called \textit{$K$-functional}, which is defined as
\begin{equation}\label{rel_Kfunctional}
K_{r,\varphi}(f,t^r)_p = \inf_h\{\|f-h\|_p+t^r\|\varphi^rh^{(r)}\|_p\,:\,h^{(r-1)}\in A.C._\text{loc}\}\;,
\end{equation}
where $h^{(r-1)}\in A.C._\text{loc}$ means that $h$ is $r-1$ times differentiable and $h^{(r-1)}$ is absolutely continuous on every closed finite interval. 
If $f\in L^p([a,b])$, then
\begin{equation}\label{rel_Kfunctional_equiv}
M^{-1}\omega_{\varphi}^r(f,t)_{p} \leq  K_{r,\varphi}(f,t^r)_p \leq M\omega_{\varphi}^r(f,t)_{p}\;,\quad\text{ for } t\leq t_0\;,
\end{equation}
for some constants $M$ and $t_0$, see Theorem 6.1.1. in \citet{ditzian87}.

\subsection{Equivalence Result}
We show that the classes $\mathcal{C}(\alpha,C)$ are equivalent to classes
defined using weighted moduli of smoothness. This, in turn, will relate them to
Sobolev and Hölder classes. To make this precise, we define for constants
$\alpha>0$ and $C>0$ the following class of functions in $\H=L^2([a,b])$
\begin{equation}\label{eq:CtildeDef}
\mathcal{\tilde C}(\alpha,C) = \{f\in \H: \|f\|_{\H}\leq C\text{ and }
\omega_{\varphi}^r(f,t)_2 \leq Ct^\alpha\text{ for all } t>0\}\;,
\end{equation}
where $\varphi(x) = \sqrt{(x-a)(b-x)}$ and $r= [\alpha]+1$.

The following theorem states the equivalence of the classes
$\mathcal{C}(\alpha,C)$ and $\mathcal{\tilde C}(\alpha,C)$. This result is
an extension of Proposition 7 in \citet{roueff05} to the case where the
subspaces $V_m$ correspond to transformed polynomial classes through an
isometry $T$ which includes both a multiplication and a composition with smooth
functions. 

\begin{theo}\label{theo_approxclass}
Let $\alpha>0$. 
Suppose that Assumption \ref{ass_polynom_L2ab} holds with a linear isometry
$T:\H=L^2([a,b])\to\H'=L^2([a',b'])$ given by $Tg=\sigma\times g\circ\tau$, where
$\sigma$ is non-negative and  $[\alpha]+1$ times continuously
differentiable, and $\tau$ is $[\alpha]+1$ times continuously
differentiable with a non-vanishing first derivative. 
Then for any positive number $\alpha$, there exist positive constants $C_1$ and $C_2$ such that for all $C>0$
\begin{equation}\label{equiv_smooth_cl}
\mathcal{C}(\alpha,C_1C)\subset \mathcal{\tilde C}(\alpha,C) \subset \mathcal{C}(\alpha,C_2C)\;.
\end{equation}
where $\mathcal{\tilde C}(\alpha,C)$ is defined in~(\ref{eq:CtildeDef})
and $\mathcal{C}(\alpha,C')$ is defined in~(\ref{def_fnspace_Cuc}) with
approximation classes $(V_m)$ given by~(\ref{eq:defVm}).
\end{theo}

For short, we write $\mathcal{C}(\alpha,\cdot)\hookrightarrow \mathcal{\tilde
  C}(\alpha,\cdot)$ when there exists  $C_1>0$ such that the first
inclusion in (\ref{equiv_smooth_cl}) holds for all $C>0$. The validity of both inclusions
is denoted by the equivalence $\mathcal{C}(\alpha,\cdot)\asymp \mathcal{\tilde C}(\alpha,\cdot)$\;.

\begin{proof}[Proof of Theorem \ref{theo_approxclass}]
Weighted moduli of smoothness are used to characterize the rate of polynomial
approximations. We start by relating $\mathcal{C}(\alpha,C)$ to classes defined
by the  rate of polynomial approximations, namely
\begin{equation*}
\mathcal{\bar C}(\alpha,C) = 
\{g\in\H': \|g\|_{H'}\leq C \text{ and } \inf_{p\in\mathcal{P}_{m-1}}
\|g-p\|_{H'}\leq C m^{-\alpha}, \text{ for all } m\geq1\}\;,
\end{equation*}
where $\mathcal{P}_m$ is the set of polynomials of degree at most $m$. Indeed,
we see that, since $T$ is a linear isometry, 
\begin{align*}
&\mathcal{C}(\alpha,C) 
= \{f\in\H : \|f\|_H\leq C \text{ and } \|P_{V_m}f-f\|_H\leq C m^{-\alpha} \text{ for all } m\geq 1\}\\
&\quad= \{T^{-1}g : g\in\H', \|g\|_{H'}\leq C \text{ and } \|P_{TV_m}g-g\|_{H'}\leq
C m^{-\alpha} \text{ for all } m\geq 1\}\\
&\quad= T^{-1}\mathcal{\bar C}(\alpha,C)\;. 
\end{align*}
As stated in Corollary 7.25 in \citet{ditzian87}, we have the equivalence $\mathcal{\bar
  C}(\alpha,\cdot) \asymp \mathcal{\tilde C}'(\alpha,\cdot)$, where
$\mathcal{\tilde C}'(\alpha,C)$ is defined as $\mathcal{\tilde C}(\alpha,C)$
but with $a'$ and $b'$ replacing $a$ and $b$. Hence, it only remains to show that
\begin{equation}\label{eq:equiv_ismetry_Ctilde}
T^{-1}\mathcal{\tilde C}'(\alpha,\cdot) \asymp \mathcal{\tilde C}(\alpha,\cdot)\;.
\end{equation}
To show this, we use the assumed particular form of $T$, that is $T(g)=\sigma
\times g\circ \tau$. Since $T$ is an isometry from $\H=L^2([a,b])$ to
$\H'=L^2([a',b'])$ and $\sigma$ is non-negative, we necessarily have that
$\tau$ is a bijection from $[a',b']$ to $[a,b]$ (whose inverse bijection is 
denoted by $\tau^{-1}$) and $\sigma=1/\sqrt{\tau'\circ\tau^{-1}}$. 
Moreover the inverse isometry writes
$T^{-1}(g)= (\sigma\circ\tau^{-1})^{-1}\times g\circ\tau^{-1}$. From the
assumptions on $\tau$ we have that $\sigma$, $(\sigma\circ\tau^{-1})^{-1}$, $\tau$ and
$\tau^{-1}$ all are $[\alpha]+1$ times continuously differentiable and the two
latter's first derivative do not vanish. 
The
equivalence~(\ref{eq:equiv_ismetry_Ctilde}) then follows by 
Lemma~\ref{lem1_Ctilde_equiv} in the appendix.
\end{proof}

\noindent\textit{Example 1 (a). Exponential Mixture (continued).}
In Example 1(a) of continuous exponential mixtures, the operator $T$ is given by (\ref{eq:T_exp_a}), that is $\sigma(t) = 1/\sqrt t$ and $\tau(t) = -\log t$ and further $\H'=L^2(\rme^{-b},\rme^{-a})$. Both $\sigma$ and $\tau$ are infinitely continuously  differentiable on $[a,b]$ if $a>0$, and thus the equivalence given in (\ref{equiv_smooth_cl}) holds.

\noindent\textit{Example 1 (b). Exponential Mixture (continued).}
For the estimator exhibited in Example 1(b) for exponential mixtures, the isometry $T$ is such that $\sigma(t) = \tau(t) = 1/t$ with $a'=1/b$ and $b'=1/a$. Hence, the conclusion of Theorem \ref{theo_approxclass} holds if $a>0$.

\noindent\textit{Example 2, 3 and 4. Gamma Shape Mixture, Scale Mixture of Beta Distributions and Exponential mixture with location parameter (continued).} In the cases of Example 2, 3 and 4, the transform $T$ is the identity and hence Theorem \ref{theo_approxclass} applies. However, this result is also obtained by Corollary 7.25 in \citet{ditzian87}.

\section{Lower Bound of the Minimax Risk}\label{sec_lowerbound}
Our goal in this section is to find a lower bound of  the minimax risk
\begin{equation*}%  \label{eq:minimax_risk_def}
\inf_{\hat f\in\mathcal{S}_n} \sup_{f\in\mathcal{C}} \pi_f^{\otimes n}\|\hat f- f\|_{\H}^2
\;,  
\end{equation*}
where $\mathcal{S}_n$ is the set of all Borel functions from $\R^n$ to $\H$,
$\mathcal{C}$ denotes a subset of densities in $\H_1$ and $\pi_f^{\otimes n}$
denotes the joint distribution of the sample $(X_1,\dots,X_n)$ under
Assumption~\ref{ass_model}. We first provide a general lower bound, which is
then used to investigate the minimax rate in the specific cases of exponential
mixtures, Gamma shape mixtures and mixtures of compactly supported scale
families.

\subsection{A General Lower Bound for Mixture Densities}
We now present a new lower bound for the minimax risk of mixture density
estimation. As in Proposition~2 in \cite{roueff05}, it relies on the mixture
structure. However, in contrast with this previous result which only applies
for mixtures of discrete distributions, we will use the following lower bound
in the case of mixtures of exponential distributions, Gamma shape
mixtures and scale mixtures of compactly supported densities.

\begin{theo}[Lower bound]\label{prop_lowbound}
  Let $f_0\in\H_1$ and $f_\ast\in\H$ with $\|f_\ast\|_{\H}\leq 1$ and $f_0\pm
  f_\ast\in\H_1$ the following lower bound holds, for any $c\in(0,1)$,
\begin{align}
  &\inf_{\hat f\in\mathcal{S}_n} \sup_{f\in\{f_0,f_0\pm f_\ast\}}
  \pi_f^{\otimes n}\|f-\hat f\|_{\H}^2  \label{rel_general_lowerbound}
  \geq c\|f_\ast\|_{\H}^2 -\frac
  c{(1-c)^2}\left( \left(1 +\int| \pi_{f_\ast}(x)|\;\zeta(\rmd x) \right)^n -
    1\right)\;,%\nonumber
\end{align}
where $\pi_f^{\otimes n}$ denotes the joint distribution of the sample
$(X_1,\dots,X_n)$ under Assumption~\ref{ass_model}.
\end{theo}

\begin{proof}
  Let $f_\ast$ be as in the Theorem. We define for a fixed $\hat f\in\mathcal{S}_n$
  and any $c\in(0,1)$ the set $A = \{\|f_0-\hat f\|_{\H}\leq \frac
  c{1-c}\}$. Then, for all $\hat f\in\mathcal{S}_n$, $\sup_{f\in\{f_0,f_0\pm
    f_\ast\}} \pi_f^{\otimes n}\|f-\hat f\|_{\H}^2$ is bounded from below by
\begin{align*}
& \frac c2 \pi_{f_0+f_\ast}^{\otimes n}\|f_0+f_\ast-\hat f\|_{\H}^2+
\frac c2 \pi_{f_0-f_\ast}^{\otimes n}\|f_0-f_\ast-\hat f\|_{\H}^2 +
(1-c)\pi_{f_0}^{\otimes n}\|f_0-\hat f\|_{\H}^2 \\
& \geq \frac c2\pi_{f_0+f_\ast}^{\otimes n}\left[\1_A\|f_0+f_\ast-\hat
  f\|_{\H}^2\right] \\
  &\qquad + \frac c2 \pi_{f_0-f_\ast}^{\otimes n}\left[\1_A\|f_0-f_\ast-\hat f\|_{\H}^2\right] +
(1-c)\pi_{f_0}^{\otimes n}\|f_0-\hat f\|_{\H}^2\;. 
\end{align*}
Note that for a function $k$ defined on $\R^n$ we have
\begin{align*}
\pi_{f_0\pm f_\ast}^{\otimes n} k
%&= \int k(x_1,\ldots,x_n) \prod_{i=1}^n \pi_{f_0\pm f_\ast}(\rmd x_i) \\
&= \int k(x_1,\ldots,x_n) \prod_{i=1}^n\left[ \pi_{f_0}(x_i) \pm
  \pi_{f_\ast}(x_i) \right] \prod_{i=1}^n\zeta(\rmd x_i) \\
&= \int k(x_1,\ldots,x_n) \sum_{I,J} \left[(\pm1)^{\#J}\prod_{j\in J}\pi_{f_*}(x_j) \prod_{i\in I} \pi_{f_0}(x_i) \right]\prod_{i=1}^n\zeta(\rmd x_i) \;,
\end{align*}
where the sum is take over all sets $I$ and $J$ such that $I\cup J =\{1,\ldots,n\}$ and $I\cap J= \emptyset$. Therefore,
\begin{align*}
&\pi_{f_0+f_\ast}^{\otimes n}\left[\1_A\|f_0+f_\ast-\hat f\|_{\H}^2\right] +
\pi_{f_0-f_\ast}^{\otimes n}\left[\1_A\|f_0-f_\ast-\hat f\|_{\H}^2\right] \\
& = \sum_{I,J} \int \prod_{i\in I} \pi_{f_0}(x_i) \prod_{j\in J}\pi_{f_\ast}(x_j)\1_A 
  \left[\|f_0+f_\ast-\hat f\|_{\H}^2 + (-1)^{\#J} \|f_0-f_\ast-\hat f\|_{\H}^2\right] \;\prod_{i=1}^n\zeta(\rmd x_i)\;.
\end{align*}
Since $\|f_*\|_{\H}\leq1$ and, on $A$, $\|f_0-\hat f\|_{\H}\leq \frac c{1-c}$, we obtain that, on $A$, $\|f_0\pm f_\ast-\hat f\|_{\H}\leq \|f_0-\hat f\|_{\H} +\|f_\ast\|_{\H}\leq \frac1{1-c}.$ This implies that the absolute value of the sum in the last display taken over all sets $I$ and $J$ such that the cardinality of set $J$ is positive, $\# J\geq1$, is lower than
\begin{align*}
& \frac 2{(1-c)^2}\sum_{I,J:\# J\geq1} \int \prod_{i\in I} \pi_{f_0}(x_i) \prod_{j\in J}| \pi_{f_\ast}(x_j)| \prod_{i=1}^n\zeta(\rmd x_i) \\
&\quad = \frac 2{(1-c)^2} \sum_{I,J:\# J\geq1} \prod_{i\in I} \int
\pi_{f_0}(x_i) 
\zeta(\rmd x_i) \prod_{j\in J}\int| \pi_{f_\ast}(x_j)|\zeta(\rmd x_j) \\  
&\quad\quad = \frac 2{(1-c)^2} \left\{ \left(1 +\int| \pi_{f_\ast}(x)|\zeta(\rmd x)  \right)^n - 1\right\}
\end{align*}
Moreover, the term with $\# J=0$ writes
\begin{align*}
&\pi_{f_0}^{\otimes n}\left( \1_A (\|f_0+f_\ast-\hat f\|_{\H}^2 + \|f_0-f_\ast-\hat f\|_{\H}^2) \right)
= 2\pi_{f_0}^{\otimes n}\left( \1_A (\|f_0-\hat f\|_{\H}^2 + \|f_\ast\|_{\H}^2) \right)\;,
\end{align*}
by the Parallelogram law. By combining these results, the minimax risk is
bounded from below by
\begin{align*}
&(1-c) \pi_{f_0}^{\otimes n}\|f_0-\hat f\|_{\H}^2 + \\
&\quad c \pi_{f_0}^{\otimes n}\left[\1_A (\|f_0-\hat f\|_{\H}^2 + \|f_\ast\|_{\H}^2) \right]-\frac c{(1-c)^2}\left[\left(1 +\int|
     \pi_{f_\ast}(x)|\zeta(\rmd x) \right)^n - 1\right]\;.
\end{align*}

% \begin{align*}
%  &(1-c) \pi_{f_0}^{\otimes n}\|f_0-\hat f\|_{\H}^2 + c \pi_{f_0}^{\otimes n}\left(\1_A (\|f_0-\hat f\|_{\H}^2 + \|f_\ast\|_{\H}^2) \right)  \\
%  &\qquad\qquad\qquad-\frac c{(1-c)^2}\left( \left(1 +\int| \pi_{f_\ast}(x)|\zeta(\rmd x) \right)^n - 1\right)
% \end{align*}
Finally we see that
\begin{align*}
&(1-c) \|f_0-\hat f\|_{\H}^2 + c\1_A\left( \|f_0-\hat f\|_{\H}^2 + \|f_\ast\|_{\H}^2 \right) \\
&\qquad\qquad= c\1_A\|f_\ast\|_{\H}^2 + ((1-c) +c\1_{A})\|f_0-\hat f\|_{\H}^2\\
&\qquad\qquad\geq c\1_A\|f_\ast\|_{\H}^2 + c\1_{A^c}\\
&\qquad\qquad\geq c\|f_\ast\|_{\H}^2\;,
\end{align*}
where we used $1\geq\|f_\ast\|_{\H}^2$. 
This yields the lower bound asserted in the theorem.
\end{proof}

\subsection{Application to Polynomial Approximation Classes}

The lower bound given in~(\ref{rel_general_lowerbound})  relies on the
choice of a function $f_\ast$ such that $f_0$ and $f_0\pm f_\ast$ are in the
smoothness class of interest.  In this subsection, we give conditions
which provide a tractable choice of $\|f_\ast\|_{\H}\leq1$ for the class
$\mathcal{C}(\alpha,C)$ defined in~(\ref{def_fnspace_Cuc}). 
Following the same lines as Theorem~1 in~\cite{roueff05}, the key idea consists in restricting our
choice using the space $V_{m}^\bot$ (the orthogonal set of $V_m$ in $\H$) and
to control separately the two terms that appear in the right hand-side
of~(\ref{rel_general_lowerbound}) within this space.

% For that purpose, we define the
% following function class for any $\alpha>0$ and any positive numbers $C$ and $K$
% \begin{equation}
%   \label{eq:ClassSup}
%   \mathcal{C}_{f_0}(K,\alpha,C) = \{f\in\H:
% \|f\|_{\infty,f_0} \leq K, \|f\|_{\H}\leq C\text{ and }\|P_{V_m}f-f\|_{\H}\leq
% C\,m^{-\alpha}\text{ for all }m\geq1\}\;,
% \end{equation}
An important constraint on $f_\ast$ is that $f_0\pm f_\ast\in\H_1$. In
particular, for controlling the sign of $f_0\pm f_\ast$, we use
the following semi-norm on $\H$,
$$
\|f\|_{\infty,f_0} = \text{ess}\sup_{t\in\Theta}\frac{|f(t)|}{f_0(t)}\;,
$$
with the convention $0/0 = 0$ and $s/0=\infty$ for $s >0$. Further, for any
subspace $V$ of $\H$, we denote 
$$
K_{\infty,f_0}(V) = \sup\{\|f\|_{\infty,f_0}: f\in V, \|f\|_{\H}=1\}\;.
$$
The following lemma will serve to optimize the term $\|f_\ast\|_{\H}$ on the
right-hand side of~(\ref{rel_general_lowerbound}). It is similar to Lemma~2
in \cite{roueff05}, so we omit its proof.
\begin{lem}\label{lem_norm_f}
  Suppose that Assumption~\ref{ass_Vm} holds. Let $f_0$ be in $\H_1$, $\alpha,C_0>0$,
  $K\leq1$ and let $\mathcal{C}(\alpha,C_0)$ be defined
  by~(\ref{def_fnspace_Cuc}) with $V_m$ given by~(\ref{eq:defVm}). Let moreover
  $w\in\H$.  Then there exists $g\in\mathcal{C}(\alpha,C_0)\cap
  V_m^\bot\cap w^\bot$ such that $\|g\|_{\infty,f_0}\leq K$ and
$$
\|g\|_{\H} = \min\left( C_0\,(m+1)^{-\alpha} , \frac K{K_{\infty,f_0}(V_{m+2}\cap V_m^\bot\cap w^\bot)}\right)\;. 
$$
\end{lem}

Under Assumption \ref{ass_polynom_L2ab}, where the orthonormal functions
$\psi_k$ are related to polynomials in some space $\H' = L^2[a',b']$, the
constant $K_{\infty,f_0}(V_{m+2}\cap V_m^\bot\cap w^\bot)$ can be bounded by
$K_{\infty,f_0}(V_{m+2})$ and then using the following lemma.

\begin{lem}\label{lem_K_infty}
  Suppose that Assumption \ref{ass_polynom_L2ab} holds. Let
  $f_0$ be in $\H_1$ and suppose that 
\begin{equation}\label{rel_inft_f_norms}
  \sup\left\{ \|f\|_{\infty,f_0}\,:\,f\in\H\text{ such that
    }\sup_{t\in[a',b']}|Tf(t)|\leq 1\right\}  <\infty \;.
\end{equation}
Then there exists a constant $C_0>0$ satisfying
$$
K_{\infty,f_0}(V_{m+2})\leq C_0m\;, \quad \text{ for all } m\geq
1\;.
$$
\end{lem}

\begin{proof}
  Note that $\{Tf:f\in V_m\}$ is the set of polynomials in $\H'$ of degree at
  most $m-1$, denoted by $\mathcal{P}_{m-1}$.  Using $\|f\|_{\H} =\|Tf\|_{\H'}$
  and denoting by $B$ the left-hand side of~(\ref{rel_inft_f_norms}), we have
\begin{align*}
  K_{\infty,f_0}(V_m)
  &= \sup\{\|f\|_{\infty,f_0}: f\in V_m, \|f\|_{\H}=1 \}\\
  &\leq  B\sup\left\{\sup_{t\in[a',b']}|Tf(t)|\,:\, f\in V_m, \|f\|_{\H}=1 \right\} \\
  &= B \sup\left\{\sup_{t\in[a',b']}|p(t)| : p\in \mathcal{P}_{m-1},
    \int_{a'}^{b'}p^2(t)\rmd t=1 \right\}
\end{align*}
By the Nikolskii inequality (see e.g. \citet{devore93}, Theorem 4.2.6), there
exists a constant $C>0$ such that the latter sup is at most $Cm$.
Hence, there exists $C_0>0$ such that $K_{\infty,f_0}(V_m)  \leq C_0m$  for all $m\geq1$.
\end{proof}

Theorem~\ref{prop_lowbound} and Lemmas~\ref{lem_norm_f} and~\ref{lem_K_infty}
yield the following result.
\begin{cor}
  \label{cor:lb}
  Let $\alpha\geq1$ and $C>(b-a)^{-1/2}$. Suppose that
  Assumption~\ref{ass_polynom_L2ab} holds with an isometry $T$ satisfying the
  assumptions of Theorem~\ref{theo_approxclass}. Let $w$ be an
  $[\alpha]+1$ times continuously differentiable function defined on $[a,b]$
  and set
  \begin{equation}
    \label{eq:cor-mn-cond}
    v_m = 
    \sup_{g\in V_m^\bot,\|g\|_{\H}\leq1}
    \int | \pi_{w\,g}(x)|  \;\zeta(\rmd x)  \;.
  \end{equation}
  Then there exists a small enough $C_\ast>0$ and $C^*>0$ such that, for any
  sequence $(m_n)$ of integers increasing to $\infty$ satisfying $v_{m_n}\leq
  C_\ast\,n^{-1}m_n^{-\alpha}$, we have 
\begin{equation}
  \label{eq:cor_lb_risk}
  \inf_{\hat f\in\mathcal{S}_n} 
  \sup_{f\in\tilde{\mathcal{C}}(\alpha,C)\cap\H_1} 
\pi_f^{\otimes n}\|\hat f- f\|_{\H}^2 \geq C^*\; m_n^{-2\alpha} (1+o(1)) \;,
\end{equation}
where $\tilde{\mathcal{C}}(\alpha,C)$ is the smoothness class defined
by~(\ref{eq:CtildeDef}).
\end{cor}
\begin{remark}
  The assumption $C>(b-a)^{-1/2}$ is necessary, otherwise
  $\tilde{\mathcal{C}}(\alpha,C)\cap \H_1$ is reduced to one density for
  $C=(b-a)^{-1/2}$ and is empty for $C<(b-a)^{-1/2}$.  To see why, observe that for
  any $f\in\H_1$, by Jensen's inequality, $\|f\|_{\H}^2=\int_a^bf^2(t)\rmd t\geq
  (b-a)^{-1}$, with equality implying that $f$ is the uniform density on
  $[a,b]$.
\end{remark}
\begin{proof}
  We apply Theorem~\ref{prop_lowbound} with $f_0$ set as the uniform density on
  $[a,b]$ and $f_\ast$ chosen as follows. For some $C_0>0$ and an integer $m$
  to be determined later, we choose $f_\ast=w g$ where $g$ is given by
  Lemma~\ref{lem_norm_f} with $K=\min(1,\sup_{t\in[a,b]}|w(t)|)$. Since $g\in
  w^\bot$ and $\|g\|_{\infty,f_0}\leq K$, we get that $f_0\pm f_\ast\in\H_1$.
  
  Now we show that $\{f_0,f_0\pm f_\ast\}\subset\tilde{\mathcal{C}}(\alpha,C)$
  for a well chosen $C_0$. We have $\|f_0\|_{\H}=(b-a)^{-1/2}$ and, since the
  symmetric differences of all order vanishes on $f_0$, we get that
  $f_0\in\tilde{\mathcal{C}}(\alpha,(b-a)^{-1/2})$.  By definition of $g$ in
  Lemma~\ref{lem_norm_f} and Lemma~\ref{lem1_Ctilde_equiv} successively, we get
  that $f_\ast\in\tilde{\mathcal{C}}(\alpha,C'_1C_0)$ for some $C'_1>0$ not
  depending on $C_0$. Choosing $C_0=(C-(b-a)^{-1/2})/C'_1$, we finally get that
  $$
  \{f_0,f_0\pm f_\ast\}\subset\tilde{\mathcal{C}}(\alpha,C)\cap\H_1\;.
  $$
  By Lemma~\ref{lem_norm_f}, $\|g\|_{\H}\to0$ as $m\to\infty$ and, since $w$ is
  bounded, it implies that $\|f_\ast\|_\H\leq1$ for $m$ large enough.  Hence
  we may apply Theorem~\ref{prop_lowbound} and, to conclude the proof, it
  remains to provide a lower bound of the right-hand side
  of~(\ref{rel_general_lowerbound}) for the above choice of $f_\ast$.  Under
  the assumptions of Theorem~\ref{theo_approxclass},
  Condition~(\ref{rel_inft_f_norms}) clearly holds. So Lemma~\ref{lem_K_infty}
  and the definition of $g$ in Lemma~\ref{lem_norm_f} give that 
  $$
  \|g\|_H\leq C'_0 m^{-\alpha}\;,
  $$
  for some constant $C'_0>0$. By definition of $v_m$ and since $g\in V_m^\bot$,
  we have
  $$
  \int| \pi_{f_\ast}(x)|\;\zeta(\rmd x)\leq \|g\|_H \,v_m \leq C'_0
  m^{-\alpha}v_m\;.
  $$
  We now apply  the lower bound given by~(\ref{rel_general_lowerbound}) with $m=m_n$
  for $(m_n)$ satisfying $v_{m_n}\leq C_\ast n^{-1}m_n^{-\alpha}$. We thus obtain
  \begin{align*}
  &\inf_{\hat f\in\mathcal{S}_n} 
  \sup_{f\in\tilde{\mathcal{C}}(\alpha,C)\cap\H_1} 
  \pi_f^{\otimes n}\|\hat f- f\|_{\H}^2 \\
  &\qquad \geq c (C'_0m_n^{-\alpha})^2
  -\frac{c}{(1-c)^2}  C_\ast C'_0 m_n^{-2\alpha}(1+o(1))\\  
  &\qquad \geq C^\ast m_n^{-2\alpha}\,(1+o(1))\;,
  \end{align*}
  where the last inequality holds for some $C^\ast>0$ provided that $C_\ast$ is
  small enough. 
\end{proof}

To apply Corollary~\ref{cor:lb},
one needs to investigate the asymptotic behavior of the sequence $(v_m)$
defined in~(\ref{eq:cor-mn-cond}). The following lemma can be used to achieve
this goal.

\begin{lem}  \label{lem:v_m_projection}
Under Assumption~\ref{ass_polynom_L2ab}, if $\pi_{\cdot}(x)\in\H$ for
all $x\in\Xset$, then $v_m$ defined in~(\ref{eq:cor-mn-cond}) satisfies
  \begin{equation}
    \label{eq:v_m_projection}
v_m  \leq \int
    \|T[w\pi_{\cdot}(x)]-P_{\mathcal{P}_{m-1}}(T[w\pi_{\cdot}(x)])\|_{\H'}\;\zeta(\rmd x)  \;,
  \end{equation}
where $\mathcal{P}_{m-1}$ is the set of polynomials of degree at
most $m-1$ in $\H'$ and $P_{\mathcal{P}_{m-1}}$ denotes the orthogonal
projection in $\H'$ onto $\mathcal{P}_{m-1}$.
\end{lem}
\begin{proof}
Let $g\in V_m^\bot$ such that $\|g\|_{\H}\leq1$. Then we
have, for all $x\in\R$,
$$
\pi_{wg}(x)= \langle wg,\pi_\cdot(x)\rangle_{\H} =
\langle g, w \pi_\cdot(x)\rangle_{\H} = 
\langle Tg,T[w\pi_\cdot(x)]\rangle_{\H'}\;.
$$
Recall that $TV_m = \mathcal{P}_{m-1}$ is the set of polynomials of degree at
most $m-1$ in $\H'$.  Hence, $Tg$ is orthogonal to
$\mathcal{P}_{m-1}$, and for any $p\in\mathcal{P}_{m-1}$, we get, for all $x\in\R$,
\begin{equation}
  \label{eq:projection_trick}
 |\pi_{wg}(x)|=
|\langle Tw,T[w\pi_\cdot(x)]-p\rangle_{\H'}|\leq\|T[w\pi_\cdot(x)]-p\|_{\H'}\;,  
\end{equation}
where we used the Cauchy-Schwarz inequality and $\|Tg\|_{\H'}= \|g\|_{\H}\leq1$.
Now the bound given by~(\ref{eq:v_m_projection}) is obtained by taking $p$ equal to the projection
of $[w\pi_{\cdot}(x)]$
onto $\mathcal{P}_{m-1}$ (observe that the right-hand side
of~(\ref{eq:projection_trick}) is then minimal).
\end{proof}

\subsection{Minimax Rate for Exponential Mixtures}
In this section, we show that in the case of exponential mixtures the
orthogonal series estimator of Example 1(a) achieves the minimax rate. 

\begin{theo}\label{theo_minimax_expmix}
  Consider the exponential case, that is, let Assumption~\ref{ass_model} hold
  with $\zeta$ defined as the Lebesgue measure on $\R_+$,
  $\Theta=[a,b]\subset(0,\infty)$ and $\pi_t(x) = t\rme^{-tx}$. Let
  $C>(b-a)^{-1/2}$ and $\alpha>1$ and define $\tilde{\mathcal{C}}(\alpha,C)$ as
  in~(\ref{eq:CtildeDef}).  Then there exists $C^*>0$ such that
\begin{align}\label{rel_lowbound_expmix}
\inf_{\hat f\in\mathcal{S}_n} \sup_{f\in\tilde{\mathcal{C}}(\alpha,C)\cap \H_1}
\pi_f^{\otimes n}\|\hat f-f\|_{\H}^2 
\geq  C^*(\log n)^{-2\alpha}(1+o(1))\;.
\end{align}
\end{theo}
\begin{proof}
  Let $g_k(x)=\1\left\{x>k-\frac12\right\}$, for $k\geq1$. Then
  Assumption~\ref{ass_polynom_L2ab} holds with $\varphi_k$ and
  $T$ defined by~(\ref{eq:varphi_exp_a}) and~(\ref{eq:T_exp_a}), respectively.
  Since $a>0$, $T$ satisfies the
  assumptions of
  Theorem~\ref{theo_approxclass}.  We may thus apply Corollary~\ref{cor:lb}
  with $w=\1_{[a,b]}$.  
  Hence the minimax lower
  bound given in~(\ref{rel_lowbound_expmix}) thus follows from~(\ref{eq:cor_lb_risk}),
  provided that we have  for some constant $C'>0$, setting $m_n=C'\log n$,
\begin{equation}
  \label{eq:last_fact_to_prove}
v_{m_n}=o(n^{-1}m_n^{-\alpha})\quad\text{as $n\to\infty$}\;,
\end{equation}
where $v_m$ is defined by~(\ref{eq:cor-mn-cond}). 
Note that $\pi_t(x)=t\rme^{-x t}\1_{\R_+}(x)$. We apply
Lemma~\ref{lem:v_m_projection} to bound $v_m$. Using the definition of $T$
in~(\ref{eq:T_exp_a}), we have for all $x\geq0$,
$[T\pi_\cdot(x)](t)  =-\log t ~t^{x-1/2}$.
We write $x\in\R_+$ as the sum of its entire and decimal parts, $x=[x]+\langle x\rangle$,
and observe that, since $<x>-1/2\in[-1/2,1/2)$ and
$[a',b']=[\rme^{-b},\rme^{-a}]\subset(0,1)$, the expansion of
$t^{<x>-1/2}=\sum_{k\geq0}\alpha_k(x) (1-t)^k$ as a power series about $t=1$
satisfies $|\alpha_k(x)|=\prod_{j=1}^k|\langle x\rangle-1/2)-j|/k!\leq1$. Extending  $-\log
t$ about $t=1$, we thus get 
$-\log(t)t^{<x>-1/2}=\sum_{k\geq0}\beta_k(x) (1-t)^k$ with
$|\beta_k(x)|=|\sum_{l=1}^k\alpha_{k-l}/l|\leq 1+\log(k)$. 
For any $x<m$, we use this expansion to approximate 
$[T\pi_\cdot(x)](t)=-\log(t)t^{<x>-1/2}\times t^{[x]}$ by a
polynomial of degree $m$. Namely, we obtain
\begin{align*}
\sup_{t\in[a',b']}|[T\pi_\cdot(x)](t)-\sum_{k=0}^{m-[x]}
\beta_k(x) \,t^{k+[x]}|
&\leq \sum_{k>m-[x]}(1+\log(k))(b')^{k+[x]} 
\leq C_1 c^m\;,
\end{align*}
where we used the bound $1+\log(k)\leq C_1(c/b')^k$, valid for some constants
$C_1>0$ and $c\in(b',1)$ not depending on $x$.  This bound also applies to
$\|T\pi_\cdot(x)-P_{\mathcal{P}_{m-1}}(T\pi_\cdot(x))\|_{\H'}$ by definition of
the projection $P_{\mathcal{P}_{m-1}}$.  For $x\geq m$,
we simply observe that $|[T\pi_\cdot(x)](t)|\leq -\log(a'){b'}^{x-1/2}$. This
also provides an upper bound for
$\|T\pi_\cdot(x)-P_{\mathcal{P}_{m-1}}(T\pi_\cdot(x))\|_{\H'}$. Finally,
integrating on $x\geq0$ we get
$$
\int_{\R_+}  \|T\pi_\cdot(x)-P_{\mathcal{P}_{m-1}}(T\pi_\cdot(x))\|_{\H'}\;\rmd x 
\leq C_2\,m\,c^{m}\;,
$$
with constants $C_2>0$ and $c<1$ not depending on $m$, and this upper bound
applies to $v_m$ by Lemma~\ref{lem:v_m_projection}. This shows
that~(\ref{eq:last_fact_to_prove}) holds provided that $C'>0$ is taken small
enough. This completes  the proof.
\end{proof}

\subsection{Minimax Rate for Gamma Shape Mixtures}
In this section, we show that in the case of Gamma shape mixtures the
orthogonal series estimator of Example~4 achieves the minimax rate up to the
$\log\log n$ multiplicative term.  

\begin{theo}\label{theo_minimax_gammamix}
  Consider the Gamma shape mixture case, that is, let
  Assumption~\ref{ass_model} hold with $\zeta$ defined as the Lebesgue measure
  on $\R_+$, $\Theta=[a,b]\subset(0,\infty)$ and $\pi_t(x) =
  x^{t-1}\rme^{-x}/\Gamma(t)$. Let $C>(b-a)^{-1/2}$ and $\alpha>1$ and define
  $\tilde{\mathcal{C}}(\alpha,C)$ as in~(\ref{eq:CtildeDef}).  Then there
  exists $C^*>0$ such that
\begin{align}\label{rel_lowbound_gammamix}
\inf_{\hat f\in\mathcal{S}_n} \sup_{f\in\tilde{\mathcal{C}}(\alpha,C)\cap \H_1}
\pi_f^{\otimes n}\|\hat f-f\|_{\H}^2 
\geq  C^*(\log n)^{-2\alpha}(1+o(1))\;.
\end{align}
\end{theo}
\begin{proof}
  We proceed as in the proof of Theorem~\ref{theo_minimax_expmix}. This time we
  set $g_k(x)=\sum_{l=1}^k\tilde{c}_{k,l}t^{l-1}$ with coefficients
  $(\tilde{c}_{k,l})$ defined in
  Lemma~\ref{lem:rec_gamma_pol}. Assumption~\ref{ass_polynom_L2ab} then holds
  with $\H'=\H$ and $T$ defined as the identity operator. Applying
  Corollary~\ref{cor:lb} with $w(t)=\Gamma(t)$, we obtain the
  lower bound given in~(\ref{rel_lowbound_gammamix}) provided that 
  Condition~(\ref{eq:last_fact_to_prove}) holds with $m_n=C'\log n/\log\log n$
  for some $C'>0$. Again we use Lemma~\ref{lem:v_m_projection} to check this
  condition in the present case. To this end we must, for each $x>0$, provide a
  polynomial approximation of $w(t)\pi_t(x)= x^{t-1}\rme^{-x}$ as a function
  of $t$. Expanding the exponential function as a power series, we get
$$
\sup_{t\in[a,b]}\left|w(t)\pi_t(x)-\rme^{-x}\sum_{k=0}^{m-1}\frac{\log^k(x)}{k!}(t-1)^k\right|\leq
\rme^{-x}\sum_{k\geq m}\frac{|\log(x)|^k}{k!}c^k\;,
$$
where $c=\max(|a-1|,|b-1|)$. Let $(x_m)$ be a sequence of real numbers tending
to infinity. The right hand side of the previous display is less than
$\rme^{c|\log(x)|-x}(c|\log(x)|)^m/m!$. We use this for bounding
$\|w\pi_{\cdot}(x)-P_{\mathcal{P}_{m-1}}(w\pi_{\cdot}(x))\|_{\H}$ (recall that
$T$ is the identity and $\H'=\H$) when $x\in[\rme^{-x_m},x_m]$.  When
$x\in(0,\rme^{-x_m})$ we use that the latter is bounded by $O(1)$ and when
$x>x_m$ by $O(\rme^{-x/2}$). Hence Lemma~\ref{lem:v_m_projection} gives
that
$$
v_m=O(\rme^{-x_m})+\frac{c^m}{m!}
\int_{\rme^{-x_m}}^{x_m}\rme^{c|\log(x)|-x}\,|\log(x)|^m\,\rmd x
+O\left(\rme^{-x_m/2}\right)\;.
$$
Now observe that, as $x_m\to\infty$, separating the integral
$\int_{\rme^{-x_m}}^{x_m}$ as $\int_{\rme^{-x_m}}^1+\int_{1}^{x_m}$, we get 
$$
\int_{\rme^{-x_m}}^{x_m}\rme^{c|\log(x)|-x}\,|\log(x)|^m\,\rmd x
=O(\rme^{cx_m}x_m^m)+O(\log^m(x_m))\;.
$$
Set $x_m=c_0m$. By Stirling's formula, for $c_0>0$ small enough, we get
$v_m=O(c_1^m)$ with $c_1\in(0,1)$. We conclude as in the proof of Theorem~\ref{theo_minimax_expmix}.
% Observe that $\Gamma$ is holomorphic on the half complex plane
%   $(0,\infty)+\rmi\R$ and does not vanish on this domain. Thus $1/\Gamma$ also
%   is holomorphic on this domain. It follows that it can be expanded as
%   $\sum_k\beta_k(t-\mu)^k$ on $t\in[a,b]\subset(0,\infty)$ with $\mu=(a+b)/2$
%   and $\beta_k=O(\delta^k)$ with some $\delta>(b-a)/2$. Hence SUFFIT PAS!!!
%   \begin{equation}
%     \label{eq:1surGamma}
%     \left|\frac1{\Gamma(t)}-\sum_{k=0}^p\beta_k(t-\mu^k)\right|\leq C_1 \eta^p \;.
%   \end{equation}
\end{proof}
\subsection{Lower Bound for Compactly Supported Scale Families}

We derived in Corollary~\ref{cor:comp-support} an upper bound of the minimax
rate for estimating $f$ in $\mathcal{C}(\alpha,C)$. It is thus legitimate to
investigate whether, as in the exponential mixture case, this upper bound is
sharp for mixtures of compactly supported scale families. A direct application of Corollary~\ref{cor:lb} provides the following
lower bound, which, unfortunately, is far from providing a complete and
definite answer.

\begin{theo}\label{theo:lb_comp_supported}
  Consider the case of scale mixtures of a compactly supported density on
  $\R_+$, that is, suppose that the assumptions of
  Lemma~\ref{lem:comp-support} hold. Suppose moreover that $\pi_1$ has a 
  $k$-th derivative bounded on $\R_+$. Let $C>(b-a)^{-1/2}$ and
  $\alpha\geq1$, and define $\tilde{\mathcal{C}}(\alpha,C)$ as
  in~(\ref{eq:CtildeDef}).  Then if $k>\alpha$, 
\begin{align}\label{rel_lowbound_comp_support}
\inf_{\hat f\in\mathcal{S}_n} \sup_{f\in\tilde{\mathcal{C}}(\alpha,C)\cap \H_1}
\pi_f^{\otimes n}\|\hat f-f\|_{\H}^2 
\geq  n^{-2\alpha/(k-\alpha)}(1+o(1))\;.
\end{align}
\end{theo}
\begin{proof}
  We proceed as in the proof of Theorem~\ref{theo_minimax_expmix}, that is, we
  observe that Assumption~\ref{ass_polynom_L2ab} holds with the same choice of
  $(g_k)$ as in Lemma~\ref{lem:comp-support} and apply
  Corollary~\ref{cor:lb} with $w=\1_{[a,b]}$. Here, the lower
  bound given in~(\ref{rel_lowbound_comp_support}) is obtained by showing that
  \begin{equation}
    \label{eq:last_fact_to_prove2}
    v_{m_n}=O(n^{-1}m_n^{-\alpha})\quad\text{as $n\to\infty$}\;,
  \end{equation}
  holds with $m_n=n^{1/(k-\alpha)}$ and with $(v_m)$ defined
  by~(\ref{eq:cor-mn-cond}). Again we use Lemma~\ref{lem:v_m_projection} to
  bound $v_m$. Here $T$ is the identity operator on $\H=\H'$ and
  $\pi_t(x)=t^{-1}\pi_1(x/t)$. Let $M>0$ such that the support of $\pi_1$ is
  included in $[0,M]$. Then for $t\in[a,b]$ and $x>Mb$, $\pi_t(x)=0$. Hence
  \begin{equation}
    \label{eq:v_m_lb_comp_support1}
   \|\pi_\cdot(x)-P_{\mathcal{P}_{m-1}}(\pi_\cdot(x))\|_{\H} =0\quad\text{for
     all $x>Mb$}\;.
  \end{equation}
  We now consider the case $x\leq Mb$. 
  By the assumption on $\pi_1$ and $a$, we have that
  $t\mapsto\pi_t(x)=t^{-1}\pi_1(x/t)$ is $k$-times differentiable on
  $[a,b]$. Moreover its $k$-th derivative is bounded by $C_kx^k$ on $[a,b]$,
  where $C_k>0$ does not depend on $x$. It follows that, for any $h>0$ and
  $t\in[a,b]$, 
  $$
  \left|\Delta_{h}^k(\pi_{\cdot}(x),t)\right|\leq C_k\,k!\,(xh)^k\;,
  $$
  where $\Delta_{h}^k$ is the $k$-th order symmetric difference operator defined
  by~(\ref{eq:diff_order_r}). Observing moreover that
  $$
  \|\pi_{\cdot}(x)\|_{\H}^2=\int_a^bt^{-2}\pi_1^2(x/t)\,\rmd t\leq C'\,,
  $$
  for some $C'>0$ not depending on $x$, we get that
  $\pi_{\cdot}(x)\in\tilde{\mathcal{C}}(k,C'\vee C_k\,k!\,x^k)$.  Using
  Corollary 7.25 in \citet{ditzian87}, we thus have for a constant $C''>0$ not depending on $x$,
  \begin{equation}
    \label{eq:v_m_lb_comp_support2}
   \|\pi_\cdot(x)-P_{\mathcal{P}_{m-1}}(\pi_\cdot(x))\|_{\H} \leq
   C''(1+x^k)m^{-k} \quad\text{for all $x\leq Mb$}\;.
  \end{equation}
  Applying Lemma~\ref{lem:v_m_projection} with~(\ref{eq:v_m_lb_comp_support1})
  and~(\ref{eq:v_m_lb_comp_support2}), we obtain $v_m=O(m^{-k})$. We conclude
  that~(\ref{eq:last_fact_to_prove2}) holds with $m_n=n^{1/(k-\alpha)}$, which
  completes the proof.  
\end{proof}

Theorem~\ref{theo:lb_comp_supported} provides polynomial lower bounds of the
minimax MISE rate whereas Corollary~\ref{cor:comp-support} gives logarithmic
upper bounds in the same smoothness spaces. Hence the question of the minimax
rate is left completely open in this case. Moreover the lower bound relies on
smoothness conditions on $\pi_1$ which rule out Example 3 (for which
$\pi_1$ is discontinuous). On the other hand, the case of scale families can be
related with the deconvolution problem that has received a considerable
attention in a series of papers of the 1990's (see
e.g. \cite{zhang90,fan1991,fan1991b,fan1993,pensky:vidakovic:1999}).
The following section sheds a light on this relationship.

\subsection{Scale Families and Deconvolution}

The following lower bound is obtained from classical lower
bounds in the deconvolution problem, derived in \cite{fan1993}.

\begin{theo}\label{theo:deconv}
  Consider the case of scale mixtures on $\R_+$, that is, suppose that
  Assumption~\ref{ass_model} with $\zeta$ equal to the Lebesgue measure on
  $\R_+$, $\Theta=[a,b]\subset(0,\infty)$ and
  $\pi_t(x)=t^{-1}\pi_1(x/t)$. Denote by $\phi$ the characteristic function of
  the density $\rme^t\pi_1(\rme^t)$ on $\R$, 
  $$
  \phi(\xi)=\int \rme^{t+\rmi\xi t}\pi_1(\rme^t)\,\rmd t\;.
  $$
  Define $\tilde{T}$ as the operator $\tilde{T}(g)=\tilde{g}$, where
  $g:\R\to\R$ and $\tilde{g}(t)=t^{-1}g(\log(t))$ for all $t\in(0,\infty)$. Let
  $C>0$ and $\alpha>0$, and define $\mathcal{L}(\alpha,C)$ as the set
  containing all densities $g$ on $\R$ such that
  \begin{equation*}%    \label{eq:Lipschitz_smoothness_class}    
\left|g^{(r)}(t)-g^{(r)}(u)\right|\leq
    C|t-u|^{\alpha-r} \text{ for all $t,u\in\R$}\;,
  \end{equation*}
  where $r=[\alpha]$.   
  \begin{enumerate}[(a)]
  \item Assume that $\phi^{(j)}(t)=O(|t|^{-\beta-j})$ as $|t|\to\infty$ for
    $j=0,1,2$, where $\phi^{(j)}$ is the $j$-th derivative of $\phi$. Then
    there exists $C^*>0$ such that
\begin{align}\label{rel_lowbound_comp_support_fan-a}
\inf_{\hat f\in\mathcal{S}_n} \sup_{f\in\tilde{T}(\mathcal{L}(\alpha,C))}
\pi_f^{\otimes n}\|\hat f-f\|_{\H}^2 
\geq  C^*\;n^{-2\alpha/(2(\alpha+\beta)+1)}(1+o(1))\;.
\end{align}
  \item Assume that $\phi(t)=O(|t|^{\beta_1}\rme^{-|t|^\beta/\gamma})$ as
    $|t|\to\infty$ for some
      $\beta,\gamma>0$ and $\beta_1$, and that
      $\pi_1(u)=o(u^{-1}|\log(u)|^{-a})$ as $u\to0,\infty$ for some $a>1$. Then
    there exists $C^*>0$ such that
\begin{align}\label{rel_lowbound_comp_support_fan-b}
\inf_{\hat f\in\mathcal{S}_n} \sup_{f\in\tilde{T}(\mathcal{L}(\alpha,C))}
\pi_f^{\otimes n}\|\hat f-f\|_{\H}^2 
\geq  C^*\;\log(n)^{-2\alpha/\beta}(1+o(1))\;.
\end{align}
  \end{enumerate}
\end{theo}
\begin{proof}
  In the scale mixture case the observation $X$ can be represented as $X=\theta
  Y$, where $Y$ and $\theta$ are independent variables having density $\pi_1$ and (unknown) density $f$, respectively. By
  taking the log of the observations, the problem of estimating the density of
  $\log(\theta)$, that is $f^*(t)=\rme^t f(\rme^t)$, is a deconvolution
  problem. Hence we may apply Theorem~2 in \cite{fan1993} to obtain lower
  bounds on the nonparametric estimation of $f^*$ from
  $\log(X_1),\dots,\log(X_n)$ under appropriate assumptions on $\phi$, which is
  the characteristic function of $\log(Y)$. Let $a'=\log(a)$ and $b'=\log(b)$.
  The lower bounds in (a) and (b) above are those appearing in (a) and (b) in
  Theorem~2 of \cite{fan1993} of the minimax quadratic risk in
  $\H'=L^2([a',b'])$ for estimating $f^*$ in the Lipschitz smoothness class
  $\mathcal{L}(\alpha,C)$. Observe that $\tilde{T}$ is defined for
  all function $g:[a',b']\to\R$ by $\tilde{T}(g)=\tilde{g}$ with $\tilde{g}$
  defined on $[a,b]$ by $\tilde{g}(t)=t^{-1}g(\log(t))$, so that
  $\tilde{T}(f^*)=f$. Observing that $\tilde{T}$ is a linear operator and that
  for any $g\in\H'$, $\|\tilde{T}(g)\|_{\H}\asymp\|g\|_{\H'}$, we obtain the
  lower bounds given in~(\ref{rel_lowbound_comp_support_fan-a})
  and~(\ref{rel_lowbound_comp_support_fan-b}).
\end{proof}

As in Theorem~\ref{theo:lb_comp_supported}, the smoother $\pi_1$ is assumed,
the slower the lower bound of the minimax rate. However the lower bounds
obtained in Theorem~\ref{theo:deconv} hold for a much larger class of scale
families. Indeed, if $\pi_1$ is compactly supported, the condition induced on
$\pi_1$ in case (a) are much weaker than in
Theorem~\ref{theo:lb_comp_supported}. For instance, it holds with $\beta=k$ for
Example~3.  For an infinitely differentiable $\pi_1$ both theorems say that the
minimax rate is slower than any polynomial rate. However, in this case, case
(b) in Theorem~\ref{theo:deconv} may provide a more precise logarithmic lower
bound. It is interesting to note that, as a consequence of
\cite{beurling.malliavin.1962}, the MISE rate $(\log n)^{-2\alpha}$, which is the
rate obtained in Corollary~\ref{cor:comp-support} by the polynomial estimator
for any compactly supported $\pi_1$, is the slowest possible minimax rate
obtained in Theorem~\ref{theo:deconv}(b) for a compactly supported $\pi_1$.
Such a comparison should be regarded with care since the smoothness class in
the latter theorem is different and cannot be compared to the smoothness
classes considered in the previous results, as we explain hereafter.

The
arguments for adapting the lower bounds of Theorem~\ref{theo:deconv} also apply
for minimax upper bounds. More precisely, using the kernel estimators for the
deconvolution problem from the observations $\log(X_1)$, ..., $\log(X_n)$ and
mapping the estimator through $\tilde{T}$, one obtains an estimator of $f$
achieving the same integrated quadratic risk. The obtained rates depend on
similar assumptions on $\phi$ as those in (a) and (b), see
\cite{fan1991,fan1991b,fan1993}.  Although the scale mixture and the
convolution model are related to one another by taking the exponential (or the
logarithm in the reverse sense) of the observations, it is important to note
that, except for Theorem~\ref{theo:deconv}, our results are of different
nature.  Indeed, the upper and lower bounds in the deconvolution problem cannot
be compared with those obtained previously in the paper because there are no
possible inclusions between the smoothness classes considered in the
deconvolution problem and those defined by polynomial approximations. 

Let us examine more closely the smoothness class
$\tilde{T}(\mathcal{L}(\alpha,C))$ that appears in the lower bounds of
Theorem~\ref{theo:deconv}, inherited from the results on the deconvolution
problem. This class contains densities with non-compact supports whereas
$\tilde{\mathcal{C}}(\alpha,C)\cap\H_1$ only contains densities with supports
in $[a,b]$. Hence neither~(\ref{rel_lowbound_comp_support_fan-a})
nor~(\ref{rel_lowbound_comp_support_fan-b}) can be used for deriving minimax
rates in $\tilde{\mathcal{C}}(\alpha,C)\cap\H_1$. In fact the densities
exhibited in \cite{fan1993} to prove the lower bound have infinite support by
construction and the argument does not at all seem to be  adaptable for a class
of compactly supported densities. As for upper bounds in the deconvolution
problem, they are based on Lipschitz or Sobolev type of smoothness conditions
which are not compatible with compactly supported densities on $[a,b]$ except
for those that are smoothly decreasing close to the end points. This follows
from the fact that, in the deconvolution problem, standard estimators (kernel
or wavelet) highly rely on the Fourier behavior both of the mixing density and
of the additive noise density.  In contrast, such boundary constraints are not
necessary for densities in $\tilde{\mathcal{C}}(\alpha,C)$. For instance the
uniform density on $[a,b]$ belongs to $\tilde{\mathcal{C}}(\alpha,C)$ for all
$\alpha>0$ and $C>(b-a)^{-1/2}$, but has a Fourier transform decreasing very
slowly. A natural conclusion of this observation is that polynomial estimators
should be used preferably to standard deconvolution estimators when the mixing
density has a known compact support $[a,b]\subset(0,\infty)$. Of course this
conclusion holds for both deconvolution and scale mixture problems.

\input{support}

\input{simul}

\appendix

\section{Technical Results}

\begin{lem}\label{lem1_Ctilde_equiv}
Let $\alpha>0$, $a<b$ and $a'<b'$. Define $\mathcal{\tilde C}(\alpha,C)_{\H}$
as in~(\ref{eq:CtildeDef}) and $\mathcal{\tilde C}'(\alpha,C)_{\H}$ similarly
with $a'$ and $b'$ replacing $a$ and $b$. Let $\sigma$ be $[\alpha]+1$ differentiable on $[a,b]$ and
$\tau:[a',b']\to[a,b]$ be $[\alpha]+1$ differentiable on $[a',b']$  with a non-vanishing first derivative. Then 
$$
\left\{\sigma f\,:\,f\in\mathcal{\tilde C}(\alpha,\cdot)_{\H}\right\}\hookrightarrow\mathcal{\tilde C}(\alpha,\cdot)_{\H}
\quad\text{and}\quad
\left\{f\circ\tau\,:\,f\in\mathcal{\tilde C}(\alpha,\cdot)_{\H}\right\}\hookrightarrow\mathcal{\tilde C}'(\alpha,\cdot)_{\H}\;.
$$ 
\end{lem}

\begin{proof}
As the first embedding is the inclusion (40) in \citet{roueff05}, we only
show the second embedding.
Let $f\in\mathcal{\tilde C}(\alpha,C)$ and denote $r= [\alpha]+1$. Let
$t\in(0,1]$. By the equivalence~(\ref{rel_Kfunctional_equiv}) with the $K$-functional given in
(\ref{rel_Kfunctional}) 
there exists a function $h$ such that $h^{(r-1)} \in A.C._\text{loc}$ and
\begin{equation}
  \label{eq:Kmodulusab}
\|f-h\|_H+t^r\|\varphi^rh^{(r)}\|_{\H}
\leq 2M\omega^r_{\varphi}(f,t)_{\H} \leq 2MCt^\alpha\;,
\end{equation}
where $\varphi(x)=\sqrt{(x-a)(b-x)}$. Let us set $\tilde h=h\circ\tau$ and show
that, for some constant $K>0$ neither depending on $t$ nor $C$,
\begin{equation}
  \label{eq:Kmodulusa'b'}
\|f\circ\tau-\tilde{h}\|_{\H'}+t^r\|\tilde{\varphi}^r\;\tilde{h}^{(r)}\|_{\H'}\leq 
K\,C\,t^\alpha\;,
\end{equation}
where we defined $\tilde{\varphi}(x)=\sqrt{(x-a')(b'-x)}$, that is
the same definition as $\varphi$ with  $a'$ and $b'$ replacing $a$ and $b$.
Using again  equivalence~(\ref{rel_Kfunctional_equiv}), the
bound given in~(\ref{eq:Kmodulusa'b'}) will achieve the proof of the lemma.

Note that since $\tau'$ does not vanish, denoting  $C_1=(\inf|\tau'|)^{-1}$, for
all $g\in\H$, we have 
\begin{equation}\label{rel_norm_equiv_ftau}
\|g\circ\tau\|_{H'}\leq C_1 \|g\|_H \;.
\end{equation}
In particular we have that
\begin{equation}
  \label{eq:Kmodulusa'b'step1}
\|f\circ\tau-\tilde{h}\|_{\H'}=\|(f-h)\circ\tau\|_{\H'}\leq C_1 \|f\|_H \;.
\end{equation}
Since $\tau$ is $r$ times continuously
differentiable and $h^{(r-1)} \in A.C._\text{loc}$, we note that $\tilde
h^{(r-1)}\in A.C._\text{loc}$ with $\tilde
h^{(r)}=\sum_{j=1}^{r} \tau_j \times  h^{(j)}\circ\tau$, where the $\tau_j$'s
are continuous functions only depending on $\tau$. Hence there is a constant
$C_2>0$ only depending on $\tau$ and $r$ such that
$$
\|\tilde{\varphi}^r\;\tilde{h}^{(r)}\|_{\H'}\leq C_2\max_{j=1,\dots,r}\|\tilde{\varphi}^r\;h^{(j)}\circ\tau\|_{\H'} \;.
$$
Another simple consequence of $\tau'$ not vanishing on $[a',b']$ is that there
exists a constant $C_3>0$ such that $\tilde{\varphi}(x)\leq C_3\varphi\circ\tau(x)$ for
all $x\in[a',b']$. Using this with~(\ref{rel_norm_equiv_ftau}) in the previous
display, we get
\begin{equation}
  \label{eq:Kmodulusa'b'step2}
\|\tilde{\varphi}^r\;\tilde{h}^{(r)}\|_{\H'}\leq C_1\,C_2\,C_3\,\max_{j=1,\dots,r}\|\varphi^r\;h^{(j)}\|_{\H}\;.
\end{equation}
We shall prove that $\|\varphi^r\;h^{(j)}\|_{\H}$ appearing in the right-hand
side of the previous inequality is in fact maximized, up to
multiplicative and additive constants, at $j=r$. For $j=1,\dots,r-1$, we proceed recursively as
follows. For any $u\in (a,b)$, we have
$$
|h^{(j)}(x)| \leq \left|\int_{u}^x h^{(j+1)}(s)\rmd s\right| +
|h^{(j)}(u)|\;.
$$
Then, by Jensen's inequality,  
\begin{align*}
&\|\varphi^r\, h^{(j)}\|_{\H}\\
& \leq 
\left\{\int_{x=a}^b \varphi^{2r}(x) \left(|x-u| \;\int_{s\in[u,x]}
    \{h^{(j+1)}(s)\}^2\rmd s\right) \,\rmd x\right\}^{1/2} + \|\varphi^r\|_{\H}\;|h^{(j)}(u)|\;,
\end{align*}
where we used the convention that $[c,d]$ denotes the same segment whether $c\leq
d$ or not. By Fubini's theorem, the term between braces reads
$$
\int_{s=a}^b \{h^{(j+1)}(s)\}^2  \psi(s;u) \rmd
s\quad\text{with}\quad\psi(s;u)=\int \1_{[u,x]}(s)\;(x-a)^r(b-x)^r|x-u|  \,\rmd x\;.
$$
Let $\tilde{a}<\tilde{b}$ be two fixed numbers in $(a,b)$. It is
straightforward to show that, for some constant $C_4>0$ only depending on
$a,b,\tilde{a},\tilde{b}$, we have
$$
 \psi(s;u) \leq C_4^2 \; \varphi^{2r}(s)\quad\text{for all $u\in(\tilde{a},\tilde{b})$}\;. 
$$
The last 3 displays thus give that
$$
\|\varphi^r\, h^{(j)}\|_{\H}\leq C_4\|\varphi^r\, h^{(j+1)}\|_{\H} +  \|\varphi^r\|_{\H}\;\inf_{u\in[\tilde{a},\tilde{b}]}|h^{(j)}(u)|\;.
$$
By induction on $j$, we thus get with~(\ref{eq:Kmodulusa'b'step2}) that there
is a constant $C_5$ such that
\begin{equation}
  \label{eq:Kmodulusa'b'step3}
\|\tilde{\varphi}^r\;\tilde{h}^{(r)}\|_{\H'}\leq C_5\,\left(
  \|\varphi^r\;h^{(r)}\|_{\H} + \sum_{j=1,\dots,r-1}\inf_{u\in[\tilde{a},\tilde{b}]}|h^{(j)}(u)|\right)\;.
\end{equation}
The final step of the proof consists in bounding
$\inf_{u\in[\tilde{a},\tilde{b}]}|h^{(j)}(u)|$ for $j=1,\dots,r-1$.
Let $\delta_j=\inf_{u\in[\tilde{a},\tilde{b}]}|h^{(j)}(u)|$.
Then for any $v,v'\in[\tilde{a},\tilde{b}]$, we have
$|h^{(j-1)}(v')-h^{(j-1)}(v)|\geq \delta_j\,|v'-v|$. Suppose that $v$ is in the
first third part of the segment $[\tilde{a},\tilde{b}]$ and $v'$ in the last
third so that $|v-v'|\geq (\tilde{b}-\tilde{a})/3$. On the other hand 
$|h^{(j-1)}(v')-h^{(j-1)}(v)|\leq|h^{(j-1)}(v')|+|h^{(j-1)}(v)|$. It follows
that $|h^{(j-1)}(v')|$ and $|h^{(j-1)}(v)|$ cannot be both less than
$\delta_j\,(\tilde{b}-\tilde{a})/3$, which provides a lower bound of 
$|h^{(j-1)}|$ on at least one sub-interval of $[\tilde{a},\tilde{b}]$ of length
$(\tilde{b}-\tilde{a})/3$. Proceeding recursively we get that there exists 
a sub-interval of $[\tilde{a},\tilde{b}]$ on which $h$ is lower bounded by
$\delta_j$ multiplied by some constant. This in turns gives that 
$$
\sum_{j=1,\dots,r-1}
\inf_{u\in[\tilde{a},\tilde{b}]}|h^{(j)}(u)|\leq  C_6 \|h\|_{\H}  \;.
$$
where $C_6$ is a constant only depending on $\tilde{a},\tilde{b}$ and $r$. 
Observe that, since $f\in\tilde{C}(\alpha,C)$, we have $\|f\|_{\H}\leq C$. 
Using~(\ref{eq:Kmodulusab}),  $t\in(0,1]$ and $\|h\|_{\H}\leq \|f-h\|_{\H} +\|f\|_{\H}$ in
the last display we thus get
$$
\sum_{j=1,\dots,r-1}
\inf_{u\in[\tilde{a},\tilde{b}]}|h^{(j)}(u)|\leq  C_6 (2M+1) C \;. 
$$
Finally, this bound,~(\ref{eq:Kmodulusa'b'step3}),~(\ref{eq:Kmodulusa'b'step1})
and~(\ref{eq:Kmodulusab}) yields~(\ref{eq:Kmodulusa'b'}) and the proof 
is achieved.
\end{proof}

\begin{lem}\label{lem:rec_gamma_pol}
Let $(p_k)$ be the sequence of polynomials defined by
$p_1(t)=1$, $p_2(t)=t$, ..., $p_k(t)=t(t+1)\dots(t+k-2)$ for all $k\geq2$.
Define the coefficients $(\tilde{c}_{k,l})_{1\leq l\leq k}$ by the expansion 
formula $t^{k-1}=\sum_{l=1}^k\tilde{c}_{k,l}p_l(t)$, valid for $k=1,2,\dots$.   
Then $\tilde{c}_{1,1}=1$, and for all $k\geq2$, 
\begin{equation}
  \label{eq:rec_gamma_pol}
  \tilde{c}_{k,1}=0\,,  \tilde{c}_{k,k}=1\quad
\text{and}\quad\tilde{c}_{k,l}=\tilde{c}_{k-1,l-1}-(l-1)\tilde{c}_{k-1,l}
\quad\text{for all}\quad l=2,\dots,k-1\;.
\end{equation}
Moreover, we have, for all $k\geq1$,
\begin{equation}
  \label{eq:rec_gamma_pol_sup_bound}
\sum_{l=1}^k|\tilde{c}_{k,l}|\leq k!\;.
\end{equation}
\end{lem}
\begin{proof}
  By definition of $p_l$, we have $tp_l(t)=p_{l+1}(t)-(l-1)p_l(t)$ for any
  $l\geq1$. Hence, for any $k\geq2$, writing
  $t^{k-1}=tt^{k-2}=\sum_l\tilde{c}_{k-1,l}tp_l(t)$, we
  obtain~(\ref{eq:rec_gamma_pol}). 

  We now prove~(\ref{eq:rec_gamma_pol_sup_bound}). It is obviously true for
  $k=1$. From~(\ref{eq:rec_gamma_pol}), it follows that, for all $k\geq 2$,
$$
\sum_{l=1}^k|\tilde{c}_{k,l}|\leq \sum_{l=1}^{k-1}l\,|\tilde{c}_{k-1,l}|+1\;.
$$
Bounding $l$ inside the last sum by $(k-1)$
%  and 1 by
% $\sum_{l=1}^{k-1}l\,|\tilde{c}_{k-1,l}|$, we get
% $\sum_{l=1}^k|\tilde{c}_{k,l}|\leq k
% \sum_{l=1}^{k-1}|\tilde{c}_{k-1,l}|$. Hence, it follows the
yields~(\ref{eq:rec_gamma_pol_sup_bound}).
\end{proof}
\bibliographystyle{plainnat}
\bibliography{biblio}

\end{document}

%% file: support.tex
\section{Support Estimation}\label{sec_support}

A basic assumption of our estimation approach is that the mixing density $f$ belongs to  $\H=L^{2}[a,b]$.  However, in practice  the exact interval $[a,b]$ is generally unknown. 
To compass this problem, we propose an estimator of the support of the mixing density $f$, or more precisely of the support of $Tf$. It can be shown that the support estimator is consistent when it is based on an estimator $T\hat f_{n,m_n}$, which is a polynomial, and  $Tf$ behaves as follows on the bounds of the support interval.

Denote by  $[a_0,b_0]$ the smallest interval such that $Tf(u)=0$ for all $u\in[a',b']\backslash[a_0,b_0]$. In other words, $a_0=\inf\{u\in[a',b'],Tf(u) >0\}$ and $b_0=\sup\{u\in[a',b'], Tf(u) >0\}$. 
Furthermore, we suppose that there exist constants $D>a_0,E<b_0,D',E',\alpha'>0$  such that
\begin{align}\label{ass_Tf_support_A}
Tf(u)&\geq ((u-a_0)/D')^{\alpha'}\;, \text{ for all } u\in[a_0,D]\;,\\
Tf(u)&\geq ((b_0-u)/E')^{\alpha'}\;, \text{ for all } u\in[E,b_0]\;. 
\label{ass_Tf_support_B}
\end{align}

For fixed $\varepsilon_n, \eta_n>0$, we define the estimators $\hat a_n$ and $\hat b_n$ of the interval bounds $a_0$ and $b_0$ by
\begin{align}
\hat a_n &= \inf\left\{u\in[a',b'] :T\hat f_{n,m_n}(v)>\frac{\varepsilon_n}2 \text{ for all } v\in[u,u+\eta_n]\right\}\label{def_estim_support_a}\\
\hat b_n &= \sup\left\{u\in[a',b'] :T\hat f_{n,m_n}(v)>\frac{\varepsilon_n}2 \text{ for all } v\in[u-\eta_n,u]\right\}\;.\label{def_estim_support_b}
\end{align}
Roughly, these estimators take the smallest and largest value where the estimator $T\hat f_{n,m_n}$ exceeds $\varepsilon_n/2$, by disregarding side-effects of size $\eta_n$.  
For a convenient choice of the sequences $(\varepsilon_n)_n$ and $(\eta_n)_n$ these estimators are consistent.

\begin{prop}\label{prop_support}
Let $\hat f_{n,m_n}$ be the density estimator defined in (\ref{def_proj_estim}) under Assumption \ref{ass_polynom_L2ab}  with $\alpha>1/2$. Suppose that $f$ verifies (\ref{ass_Tf_support_A}-\ref{ass_Tf_support_B}) for appropriate  constants $D>a_0,E<b_0,D',E',\alpha'>0$.
Assume that there are sequences $m_n\to\infty$, $\varepsilon_n\to0$ and $\eta_n\to0$ such that
\begin{align*}
	  \E\left\|\hat f_{n,m_n}-f\right\|_{\H}^2  =O\left(m_n^{-2\alpha}\right)\;,\qquad
	 \varepsilon_n^{-1}=o\left(m_n^{(2\alpha-1)/(2+1/\alpha')}\right)\;,\qquad
	 \eta_n=O\left(\varepsilon_n^{1/\alpha'}m_n^{-1}\right)\;.
\end{align*}
Then the estimators $\hat a_n$ and $\hat b_n$ defined by (\ref{def_estim_support_a}) and (\ref{def_estim_support_b}) are consistent for the support bounds $a_0$ and $b_0$. More precisely, as $n\to\infty$,
	\begin{align*}
 		&(\hat a_n-a_0)_+ = O_{P}\left(\varepsilon_n^{1/\alpha'}\right)\quad\text{ and } \quad(\hat a_n-a_0)_-=O_{P}\left(\varepsilon_n^{1/\alpha'}m_n^{-1}\right)\;,\\
		&(\hat b_n-b_0)_+ = O_{P}\left(\varepsilon_n^{1/\alpha'}m_n^{-1}\right)\quad\text{ and } \quad(\hat b_n-b_0)_-=O_{P}\left(\varepsilon_n^{1/\alpha'}\right)\;.
 	\end{align*}
\end{prop}

\begin{proof}
First we consider $(\hat a_n-a_0)_+$. 
We set $\delta_n=MD'\varepsilon_n^{1/\alpha'}$ for some $M>1$ and denote
\begin{align*}
A_n&=\{(\hat a_n-a_0)_+>\delta_n\}
=\{\hat a_n> a_0+\delta_n\} \\
&=\left\{ \forall u\in [a', a_0+\delta_n] ~\exists v\in[u,u+\eta_n] \text{ such that }T\hat  f_{n,m_n}(v) \leq\frac{\varepsilon_n}2 \right\}\;.
\end{align*}
As $T\hat f_{n,m_n}$ is a  polynomial of degree $m_n$, $T\hat f_{n,m_n}$  has at most $m_n$ intersections with any constant function. Hence, 	 the number of subintervals of $[a',b']$ where $T\hat f_{n,m_n}$ exceeds $\varepsilon/2$ for any fixed $\varepsilon>0$ is bounded by $m_n$. On $A_n$, all such intervals included in $[a',a_0+\delta_n]$ are at most of size $\eta_n$. Thus, on $A_n$,
	\begin{align*}
 		 \int_{a'}^{a_0+\delta_n}\1\left\{T\hat f_{n,m_n}(u)>\frac{\varepsilon_n}2\right\}\rmd u
		 \leq m_n\eta_n\;.
 	\end{align*}
	It follows, that on $A_n$,
	\begin{align*}
 		\int_{a'}^{a_0+\delta_n}\1\left\{T\hat f_{n,m_n}(u)\leq \frac{\varepsilon_n}2\right\}\rmd u
		&\geq a_0+\delta_n-a'-m_n\eta_n\;,
 	\end{align*}
	and thus
	\begin{align*}
 		\int_{a_0+D'\varepsilon_n^{1/\alpha'}}^{a_0+\delta_n}\1\left\{T\hat f_{n,m_n}(u)\leq \frac{\varepsilon_n}2\right\}\rmd u
		&\geq \delta_n-m_n\eta_n-D'\varepsilon_n^{1/\alpha'}\;.
 	\end{align*}
	For large $n$  such that $\delta_n^{1/\alpha'}<D$ and since $Tf>\varepsilon_n$ on $[a_0+D'\varepsilon_n^{1/\alpha'},D]$ by~(\ref{ass_Tf_support_A}), we obtain on $A_n$,
	\begin{align*}
 		 &\int_{a_0+D'\varepsilon_n^{1/\alpha'}}^{a_0+\delta_n}\1\left\{T\hat f_{n,m_n}(u)\leq\frac{\varepsilon_n}2\right\}\rmd u
		 \leq \int_{a_0+D'\varepsilon_n^{1/\alpha'}}^{a_0+\delta_n}\1\left\{|T\hat f_{n,m_n}(u)-Tf(u)|>\frac{\varepsilon_n}2\right\}\rmd u \\
		 &\quad\leq \frac4{\varepsilon_n^{2}}\int_{a_0+D'\varepsilon_n^{1/\alpha'}}^{a_0+\delta_n}|T\hat f_{n,m_n}(u)-Tf(u)|^{2}\rmd u
		 \leq\frac{4}{\varepsilon_n^{2}}\|T\hat f_{n,m_n}-Tf\|_{\H'}^{2}
		 =\frac{4}{\varepsilon_n^{2}}\|\hat f_{n,m_n}-f\|_{\H}^{2}\;.		 
 	\end{align*}
	For sufficiently large $M$ we have $m_n\eta_n<\delta_n-D'\varepsilon_n^{1/\alpha'}$. Then, it follows by Markov's inequality that
	\begin{align*}
 		\P((\hat a_n- a_0)_+>\delta_n)
		&\leq \P\left( \frac{4}{\varepsilon_n^{2}}\|\hat f_{n,m_n}-f\|_{\H}^{2} \geq\delta_n-m_n\eta_n-D'\varepsilon_n^{1/\alpha'} \right) \\
		&\leq \frac{4\E[\|\hat f_{n,m_n}-f\|_{\H}^{2}]}{\varepsilon_n^{2}(\delta_n-m_n\eta_n-D'\varepsilon_n^{1/\alpha'})}
		\longrightarrow0\;,\quad n\to\infty\;,
 	\end{align*}	
	by the assumptions on $(\varepsilon_n)_n$  and $ \E\left\|\hat f_{n,m_n}-f\right\|_{\H}^2$ and as $\delta_n=MD'\varepsilon_n^{1/\alpha'}$. Thus $(\hat a_n- a_0)_+=O_{P}(\delta_n) = O_{P}\left(\varepsilon_n^{1/\alpha'}\right)$.

To investigate $(\hat a_n-a_0)_-$  put $\delta_n=M'\eta_n$ for some $M'>1$. By using that $Tf=0$ on $[a,a_0]$, we have
	\begin{align*}
 \P((\hat a_n-a_0)_->\delta_n)&=		\P(\hat a_n< a_0-\delta_n)\\
		&=\P\left( \exists x\in [a', a_0-\delta_n[ : \int_{x}^{x+\eta_n}\1\left\{T \hat f_{n,m_n}(u)>\frac{\varepsilon_n}2\right\}\rmd u =\eta_n\right) \\
		&\leq \P\left(\int_{a'}^{a_0}\1\left\{T\hat f_{n,m_n}(u)>\frac{\varepsilon_n}2\right\}\rmd u \geq\eta_n\right) \\
		&= \P\left(\int_{a'}^{a_0}\1\left\{|T\hat f_{n,m_n}(u)-Tf(u)|^{2}>\frac{\varepsilon_n^{2}}4\right\}\rmd u >\eta_n\right) \\
		&\leq \P\left(\frac4{\varepsilon_n^{2}} \int_{a}^{a_0}|T\hat f_{n,m_n}(u)-Tf(u)|^{2}\rmd u >\eta_n\right) \\
		&\leq\frac{4\E[\|\hat f_{n,m_n}-f\|^{2}_{\H}]}{\eta_n\varepsilon_n^{2}}\\
		&\longrightarrow0\;, 
	\end{align*}
	where again we applied Markov's inequality. Consequently, $(\hat a_n-a_0)_-=O_P(\eta_n)=O_{P}\left(\varepsilon_n^{1/\alpha'}m_n^{-1}\right).$
	
	By symmetry, the properties on $\hat b_n$ stated in the proposition hold as well.
	\end{proof}

By Theorem \ref{theo_mise_upperbound} the proposition applies to Example 1 (a) and 3 with $m_n=A\log n$ and to Example 1~(b) and 2 with $m_n=A\log n/\log\log n$.

%% file: simul.tex
\section{Numerical results}\label{sec_simul}

\begin{table}
 \caption{
 Estimated MISE  (and standard deviation)  of estimator $\hat f_{m,n}$ with $m=5$ in six different mixture settings when the mixing density $f$ is a Beta distribution. }\label{table_mise_betamixing}
 \begin{tabular}{lc|cccccccc}
& &\multicolumn{8}{c}{$n$}\\
&&$10^2$&$10^3$&$10^4$&$10^5$&$10^6$&$10^7$&$10^8$&$10^9$\\\hline
 Exp. (a)
&MISE&0.72&0.69&0.62&0.48&0.26&0.058&8.4e-03&1.1e-03\\
&sd &(0.16)&(0.17)&(0.16)&(0.19)&(0.18)&(0.085)&(0.10)&(1.4e-03)\\\hline
 Exp. (b)
&MISE &0.61&0.52&0.35&0.21&0.084&0.015&2.0e-03&4.4e-04\\
&sd &(0.21)&(0.25)&(0.25)&(0.21)&(0.12)&(0.027)&(2.9e-03)&(2.7e-04)\\\hline
 Gamma
&MISE &0.58&0.47&0.31&0.12&0.020&3.4e-03&1.5e-03&1.3e-03\\
&sd &(0.20)&(0.22)&(0.21)&(0.13)&(0.024)&(3.0e-03)&(4.5e-04)&(4.9e-05)\\\hline
 Uniform 
&MISE &0.32&0.10&0.015&2.9e-03&1.5e-03&1.3e-03&1.3e-03&1.3e-03\\
&sd &(0.25)&(0.12)&(0.018)&(2.1e-03)&(2.6e-04)&(4.4e-05)&(1.2e-05)&(3.6e-06)\\\hline
 Beta
&MISE &0.46&0.19&0.035&5.4e-03&1.7e-03&1.4e-03&1.3e-03&1.3e-03\\
&sd &(0.27)&(0.17)&(0.040)&(5.3e-03)&(6.5e-04)&(8.6e-05)&(2.2e-05)&(5.8e-06)\\\hline
Exp. loc.
&MISE &0.55&0.47&0.29&0.11&0.015&2.9e-03&1.5e-03&1.3e-03\\
&sd &(0.20)&(0.23)&(0.21)&(0.11)&(0.018)&(1.9e-03)&(2.6e-04)&(5.3e-05)\\\hline
 \end{tabular}
    \end{table}

\begin{figure}[t]
 \begin{center}
 \begin{tabular}{ccc}
 \includegraphics[scale=0.325]{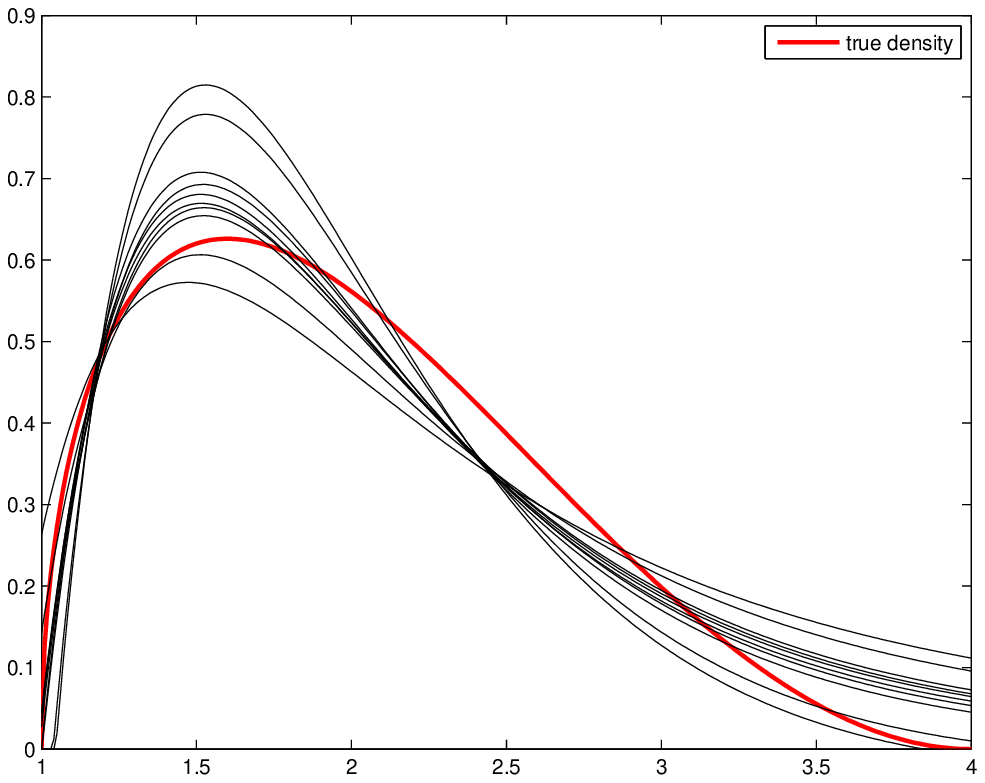}
& \includegraphics[scale=0.325]{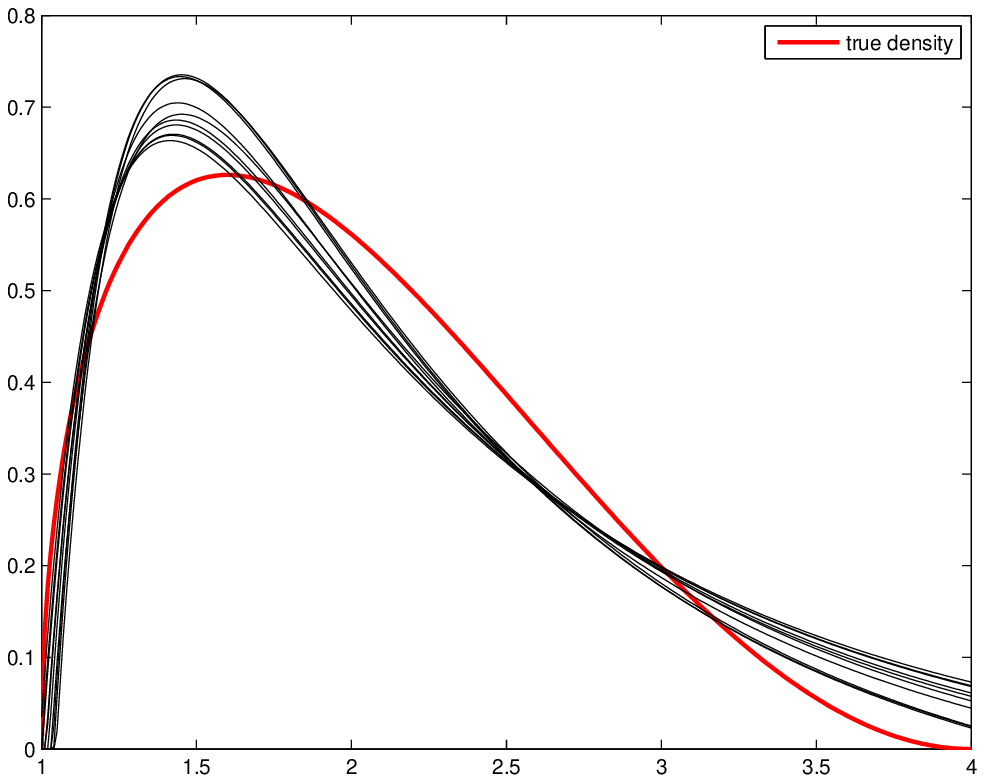}
%& \includegraphics[scale=0.4325]{fig_gamma_n105_m4.eps}
& \includegraphics[scale=0.325]{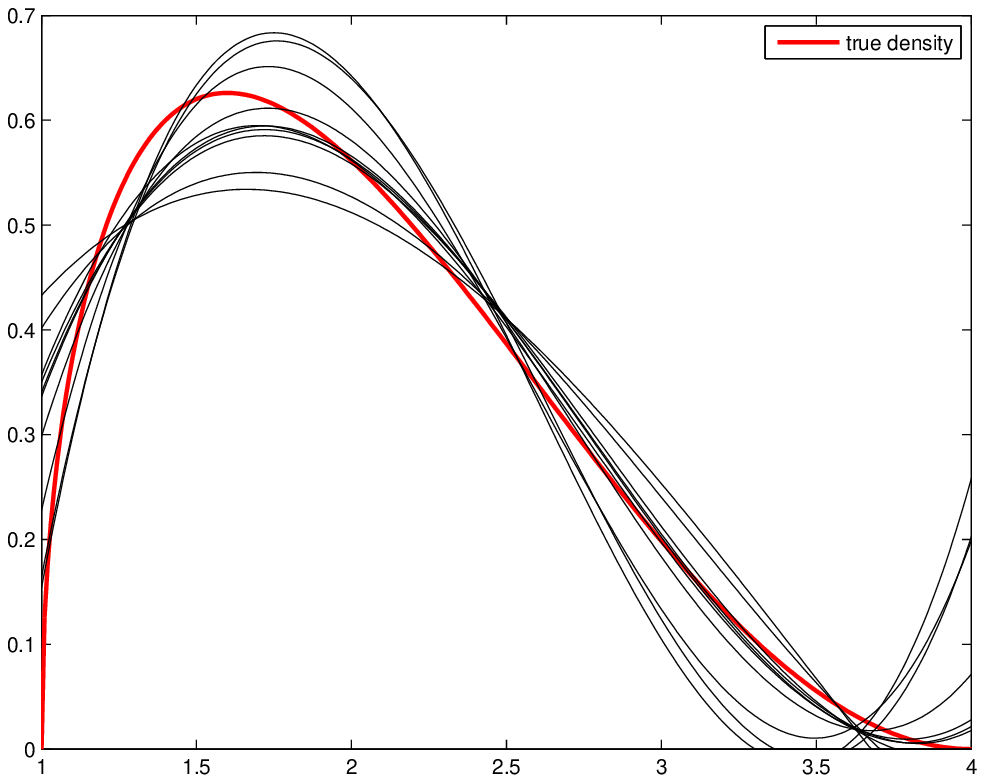}
\\
 (a) Exp. mixture (a), $m_{\text{best}}=3$&
 (b) Exp. mixture (b), $m_{\text{best}}=3$& 
 (c) Gamma mixture, $m_{\text{best}}=4$\\
 \includegraphics[scale=0.325]{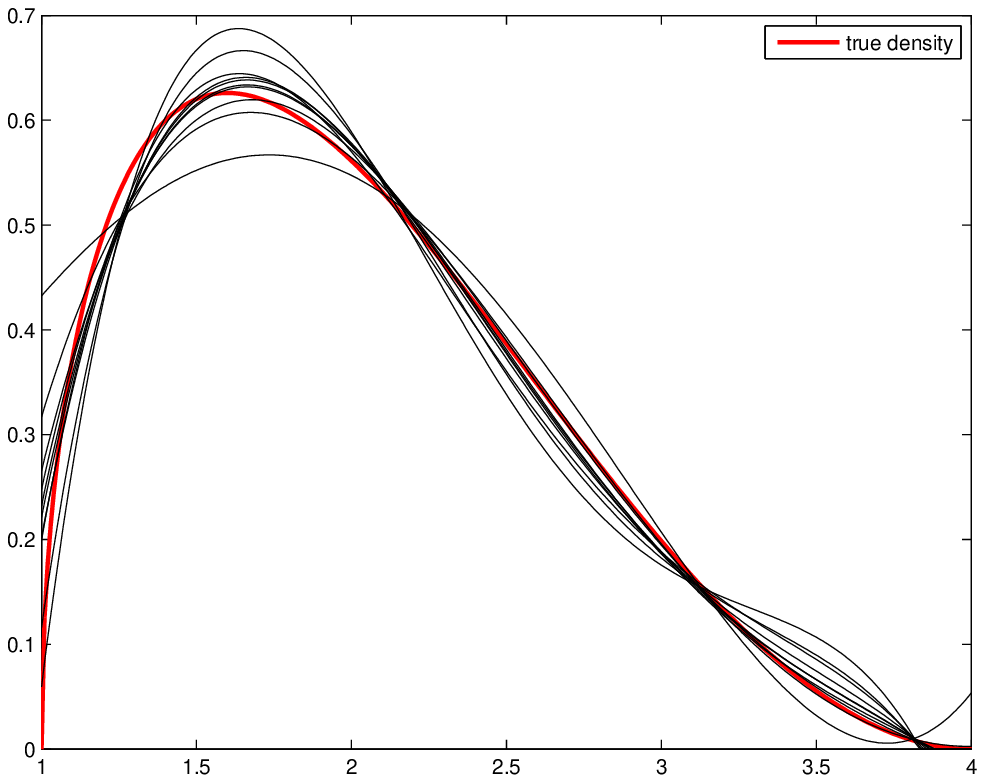}
& \includegraphics[scale=0.325]{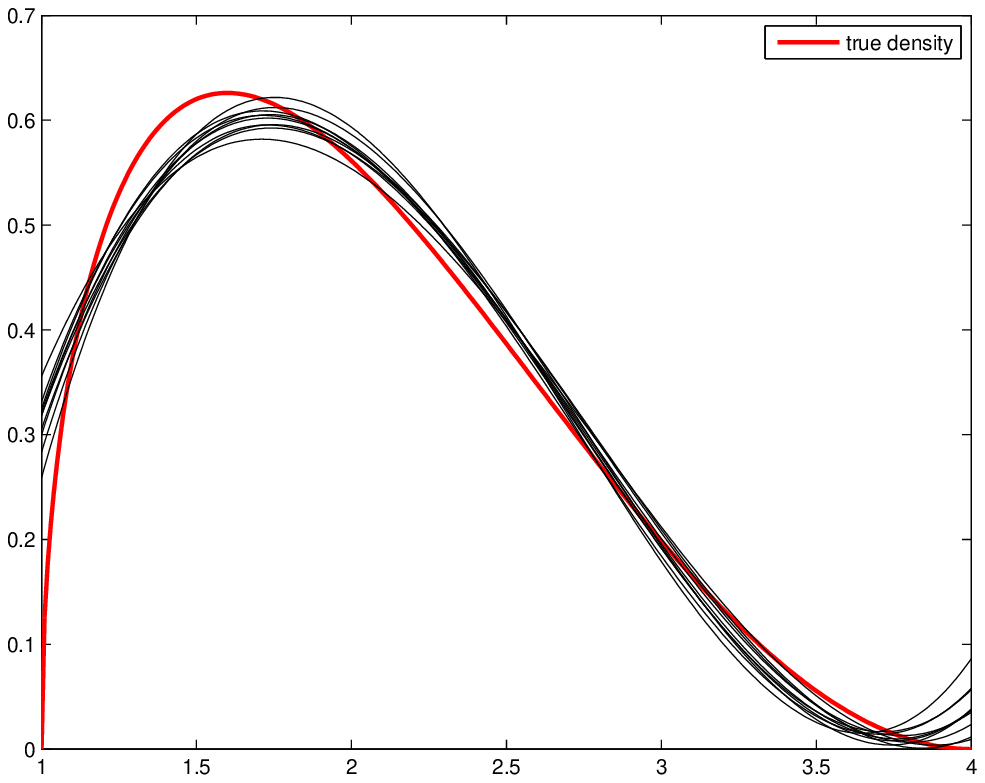}
& \includegraphics[scale=0.325]{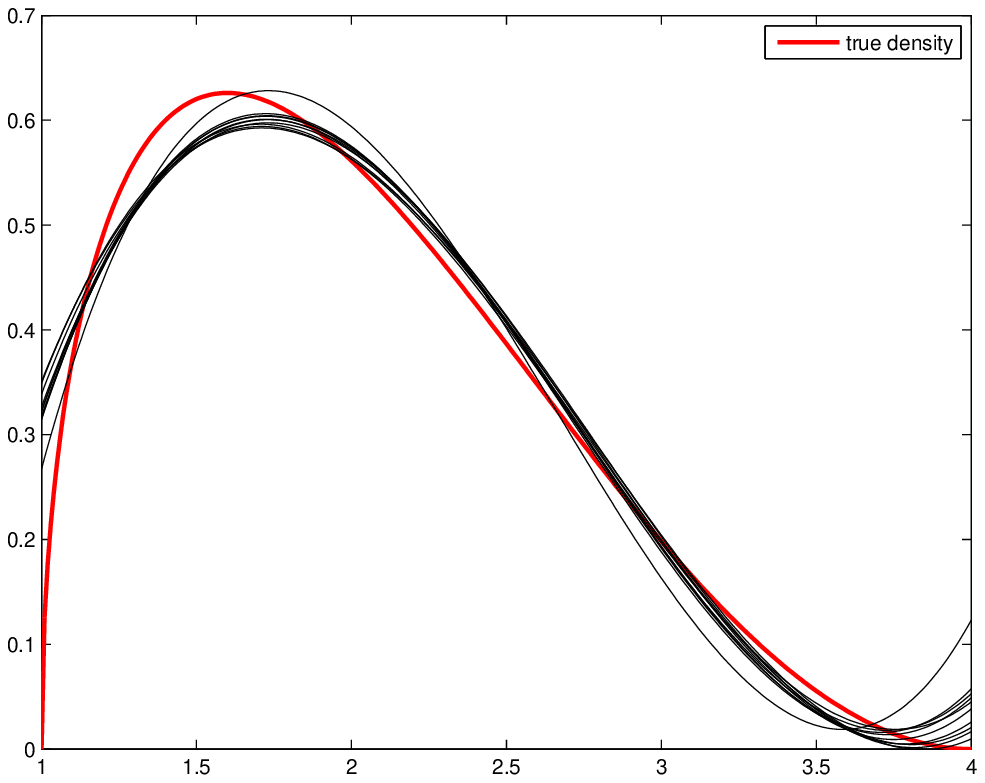}
\\
(d) Uniform mixture, $m_{\text{best}}=5$ 
&  (e) Beta mixture, $m_{\text{best}}=4$
& (f) Exp. location, $m_{\text{best}}=4$\\
  \end{tabular}
 \caption{10 estimators $\hat f_{m_{\text{best}},n}$ (black) in six different settings 
  with 
 $n=10^{5}$ 
 when the  mixing density $f$ is a Beta distribution (red).}\label{fig_estim_beta_ne5_mbest}
 \end{center}
 \end{figure}

A simulation study is conducted to evaluate the performance of the estimator on finite datasets. Six different mixture settings are considered, namely the exponential mixture from Example 1 (a) and 1 (b), the Gamma shape mixture from Example 2, the 
uniform mixture and the Beta mixture with $k=4$ from Example 3 and the exponential mixture with a location parameter from Example 4.

We consider the case where the mixing density $f$ is the Beta distribution on the interval $[1,4]$ with parameters $\alpha=3/2$ and $\beta=3$. Remark that for the exponential mixture setting of Example 1 (b) we cannot take a mixing distribution with support $[a,b]$ with $a=0$, since $b'=1/a$ must be finite. 

  For every mixture setting,  the estimator  $\hat f_{m,n}$ with $m=5$ is computed on a large number of datasets (for sample sizes $n$ varying from 100 to $10^9$) and the corresponding MISE is evaluated. Table \ref{table_mise_betamixing} gives the mean values of the different MISE and the associated standard deviations. Obviously, in all six settings the MISE decreases when $n$ increases.
Note that  in the last four settings, where the mixing density $f$ is approximated in  the same polynomial  basis,   the MISE tends to the same value, which is obviously the squared bias of the estimator when $m=5$. In the exponential mixture settings, different values are obtained because different bases are used to approach $f$.
The exponential mixture setting from Example 1 (a) always has the largest mean MISE value, while the uniform and the Beta mixtures are doing best.
  
  Figure \ref{fig_estim_beta_ne5_mbest} illustrates the estimator $\hat f_{m,n}$  when $n=10^5$ and where the order $m$ is the value  minimizing the MISE when $n=10^5$, say $m_\text{best}$. The values of $m_\text{best}$ have been obtained by extra simulations. We see that in the first two settings, we only have $m_\text{best}=3$ and the estimator seems slightly biased. On the contrary, the uniform mixture setting allows for the best approximation with $m_\text{best}=5$. 
  
%Next we consider the more involved case, where the mixing density $f$ is a
%mixture of two Gaussian distributions restricted to the interval $[1,2]$. Again, for every mixture setting a large number of datasets is  generated and we computed the estimator  $\hat f_{m,n}$ and the associated MISE for $m$ going from 1 to 6. Table \ref{table_mise_gaussmix} shows for every $n$ the order $m$ that minimizes the MISE, denoted by $m_\text{best}$, and the associated estimated MISE value and its standard deviation. Clearly, the value of the best order $m_\text{best}$ is increasing with the sample size $n$ in all six settings. Interestingly, there are jumps, that is, for the uniform mixture, for example,  $m_\text{best}$ takes the values $1,3$ and 6, but leaves out 2, 4 and 5.  
%For sample sizes up to $10^8$, the first two exponential settings achieve the best MISE values. However, for $n=10^9$, the uniform and the Beta setting are doing better, probably due to the high values of $m_{\text{best}}$. 
%
%Figure \ref{fig_estim_gauss_ne10_m6} illustrates the estimator when $n=10^{10}$ and $m=6$. Again, the estimator has the best performance in the uniform and the Beta mixture settings, where the estimator captures quite well the  form of the Gaussian mixture and the   variance of the estimators is 0. In the first two exponential settings the estimator has low variance but substantial bias.  And in the remaining cases (c) and (f)   the estimator is biased and has non-negligible variance.